
\pretolerance=500 \tolerance=1000  \brokenpenalty=5000

\catcode`\@=11


\font\eightrm=cmr8
\font\eighti=cmmi8
\font\eightsy=cmsy8
\font\eightbf=cmbx8
\font\eighttt=cmtt8
\font\eightit=cmti8
\font\eightsl=cmsl8
\font\sevenrm=cmr7
\font\seveni=cmmi7
\font\sevensy=cmsy7
\font\sevenbf=cmbx7

\font\sixrm=cmr6
\font\sixi=cmmi6
\font\sixsy=cmsy6
\font\sixbf=cmbx6

\skewchar\eighti='177 \skewchar\sixi='177
\skewchar\eightsy='60 \skewchar\sixsy='60

\def\tenpoint{%
   \textfont0=\tenrm \scriptfont0=\sevenrm
   \scriptscriptfont0=\fiverm
   \def\rm{\fam\z@\tenrm}%
   \textfont1=\teni  \scriptfont1=\seveni
   \scriptscriptfont1=\fivei
   \def\oldstyle{\fam\@ne\teni}\let\old=\oldstyle
   \textfont2=\tensy \scriptfont2=\sevensy
   \scriptscriptfont2=\fivesy
   \textfont\itfam=\tenit
   \def\it{\fam\itfam\tenit}%
   \textfont\slfam=\tensl
   \def\sl{\fam\slfam\tensl}%
   \textfont\bffam=\tenbf
   \scriptfont\bffam=\sevenbf
   \scriptscriptfont\bffam=\fivebf
   \def\bf{\fam\bffam\tenbf}%
   \textfont\ttfam=\tentt
   \def\tt{\fam\ttfam\tentt}%
   \abovedisplayskip=12pt plus 3pt minus 9pt
   \belowdisplayskip=\abovedisplayskip
   \abovedisplayshortskip=0pt plus 3pt
   \belowdisplayshortskip=4pt plus 3pt
   \smallskipamount=3pt plus 1pt minus 1pt
   \medskipamount=6pt plus 2pt minus 2pt
   \bigskipamount=12pt plus 4pt minus 4pt
   \normalbaselineskip=12pt
   \setbox\strutbox=\hbox{\vrule height8.5pt depth3.5pt width0pt}%
   \let\bigf@nt=\tenrm
   \let\smallf@nt=\sevenrm
   \normalbaselines\rm}

\def\eightpoint{%
   \textfont0=\eightrm \scriptfont0=\sixrm
   \scriptscriptfont0=\fiverm
   \def\rm{\fam\z@\eightrm}%
   \textfont1=\eighti  \scriptfont1=\sixi
   \scriptscriptfont1=\fivei
   \def\oldstyle{\fam\@ne\eighti}\let\old=\oldstyle
   \textfont2=\eightsy \scriptfont2=\sixsy
   \scriptscriptfont2=\fivesy
   \textfont\itfam=\eightit
   \def\it{\fam\itfam\eightit}%
   \textfont\slfam=\eightsl
   \def\sl{\fam\slfam\eightsl}%
   \textfont\bffam=\eightbf
   \scriptfont\bffam=\sixbf
   \scriptscriptfont\bffam=\fivebf
   \def\bf{\fam\bffam\eightbf}%
   \textfont\ttfam=\eighttt
   \def\tt{\fam\ttfam\eighttt}%
   \abovedisplayskip=9pt plus 3pt minus 9pt
   \belowdisplayskip=\abovedisplayskip
   \abovedisplayshortskip=0pt plus 3pt
   \belowdisplayshortskip=3pt plus 3pt
   \smallskipamount=2pt plus 1pt minus 1pt
   \medskipamount=4pt plus 2pt minus 1pt
   \bigskipamount=9pt plus 3pt minus 3pt
   \normalbaselineskip=9pt
   \setbox\strutbox=\hbox{\vrule height7pt depth2pt width0pt}%
   \let\bigf@nt=\eightrm
   \let\smallf@nt=\sixrm
   \normalbaselines\rm}
\tenpoint

\font\tencal=eusm10

\font\sevencal=eusm7

\font\fivecal=eusm5
\newfam\calfam
\textfont\calfam=\tencal
\scriptfont\calfam=\sevencal
\scriptscriptfont\calfam=\fivecal
\def\cal#1{{\fam\calfam\relax#1}}

\def\pc#1{\bigf@nt#1\smallf@nt}

\catcode`\;=\active
\def;{\relax\ifhmode\ifdim\lastskip>\z@\unskip\fi \kern\fontdimen2
 -1.2 \fontdimen3 \string;}

\catcode`\:=\active
\def:{\relax\ifhmode\ifdim\lastskip>\z@\unskip\fi\penalty\@M\
\fi\string:}

\catcode`\!=\active
\def!{\relax\ifhmode\ifdim\lastskip>\z@ \unskip\fi\kern\fontdimen2
 -1.1 \fontdimen3 \string!}

\catcode`\?=\active
\def?{\relax\ifhmode\ifdim\lastskip>\z@ \unskip\fi\kern\fontdimen2
 -1.1 \fontdimen3 \string?}

\frenchspacing

\def\og{\leavevmode\raise.3ex\hbox{$\scriptscriptstyle
\langle\!\langle\,$}}
\def\fg{\leavevmode\raise.3ex\hbox{$\scriptscriptstyle
\,\rangle\!\rangle$}}

\def\pointir{\unskip . --- \ignorespaces}


\def\Medbreak{\vskip-\lastskip\medbreak}

\def\rem#1\endrem{%
\Medbreak {\it#1\unskip} : }

\long\def\thm#1 #2\enonce#3\endthm{%
\Medbreak {\pc#1} {#2\unskip}\pointir{\it #3}\medskip}

\def\decale#1{\smallbreak\hskip 28pt\llap{#1}\kern 5pt}
\def\decaledecale#1{\smallbreak\hskip 34pt\llap{#1}\kern 5pt}

\let\@ldmessage=\message

\def\message#1{{\def\pc{\string\pc\space}%
\def\'{\string'}\def\`{\string`}%
\def\^{\string^}\def\"{\string"}%
\@ldmessage{#1}}}


\def\up#1{\raise 1ex\hbox{\smallf@nt#1}}

\def\diagram#1{\def\normalbaselines{\baselineskip=0pt
\lineskip=5pt}\matrix{#1}}

\def\longmapsto#1{\mapstochar\mathrel{\joinrel \kern-0.2mm\hbox to
#1mm{\rightarrowfill}}}

\catcode`\@=12

\showboxbreadth=-1  \showboxdepth=-1


\message{`lline' & `vector' macros from LaTeX}

\def\Grille{\setbox13=\vbox to 5\unitlength{\hrule width 109mm \vfill}
\setbox13=\vbox to 65mm
{\offinterlineskip\leaders\copy13\vfill\kern-1pt\hrule}
\setbox14=\hbox to 5\unitlength{\vrule height 65mm\hfill}
\setbox14=\hbox to 109mm{\leaders\copy14\hfill\kern-2mm \vrule height
65mm}
\ht14=0pt\dp14=0pt\wd14=0pt \setbox13=\vbox to 0pt
{\vss\box13\offinterlineskip\box14} \wd13=0pt\box13}

\def\rule(#1,#2)\dir(#3,#4)\long#5{%
\noalign{\leftput(#1,#2){\lline(#3,#4){#5}}}}
\def\arrow(#1,#2)\dir(#3,#4)\length#5{%
\noalign{\leftput(#1,#2){\vector(#3,#4){#5}}}}
\def\put(#1,#2)#3{\noalign{\setbox1=\hbox{%
\kern #1\unitlength \raise #2\unitlength\hbox{$#3$}}%
\ht1=0pt \wd1=0pt \dp1=0pt\box1}}

\catcode`@=11

\def\{{\relax\ifmmode\lbrace\else$\lbrace$\fi}
\def\}{\relax\ifmmode\rbrace\else$\rbrace$\fi}
\def\newcount{\alloc@0\count\countdef\insc@unt}
\def\newdimen{\alloc@1\dimen\dimendef\insc@unt}
\def\newwrite{\alloc@7\write\chardef\sixt@@n}

\newwrite\@unused
\newcount\@tempcnta
\newcount\@tempcntb
\newdimen\@tempdima
\newdimen\@tempdimb
\newbox\@tempboxa

\def\@spaces{\space\space\space\space}
\def\@whilenoop#1{}
\def\@whiledim#1\do #2{\ifdim #1\relax#2\@iwhiledim{#1\relax#2}\fi}
\def\@iwhiledim#1{\ifdim #1\let\@nextwhile=\@iwhiledim
\else\let\@nextwhile=\@whilenoop\fi\@nextwhile{#1}}
\def\@badlinearg{\@latexerr{Bad \string\line\space or \string\vector
\space argument}}
\def\@latexerr#1#2{\begingroup
\edef\@tempc{#2}\expandafter\errhelp\expandafter{\@tempc}%

\def\@eha{Your command was ignored.^^JType \space I <command> <return>
\space to replace it with another command,^^Jor \space <return> \space
to continue without it.}
\def\@ehb{You've lost some text. \space \@ehc}
\def\@ehc{Try typing \space <return> \space to proceed.^^JIf that
doesn't work, type \space X <return> \space to quit.}
\def\@ehd{You're in trouble here.  \space\@ehc}
\typeout{LaTeX error.  \space See LaTeX manual for explanation.^^J
\space\@spaces\@spaces\@spaces Type \space H <return> \space for
immediate help.}\errmessage{#1}\endgroup}
\def\typeout#1{{\let\protect\string\immediate\write\@unused{#1}}}

\font\tenln = line10
\font\tenlnw = linew10

\newdimen\@wholewidth
\newdimen\@halfwidth
\newdimen\unitlength

\unitlength =1pt

\def\thinlines{\let\@linefnt\tenln \let\@circlefnt\tencirc
\@wholewidth\fontdimen8\tenln \@halfwidth .5\@wholewidth}
\def\thicklines{\let\@linefnt\tenlnw \let\@circlefnt\tencircw
\@wholewidth\fontdimen8\tenlnw \@halfwidth .5\@wholewidth}
\def\linethickness#1{\@wholewidth #1\relax \@halfwidth .5
\@wholewidth}

\newif\if@negarg

\def\lline(#1,#2)#3{\@xarg #1\relax \@yarg #2\relax
\@linelen=#3\unitlength \ifnum\@xarg =0 \@vline \else \ifnum\@yarg =0
\@hline \else \@sline\fi \fi}
\def\@sline{\ifnum\@xarg< 0 \@negargtrue \@xarg -\@xarg \@yyarg
-\@yarg \else \@negargfalse \@yyarg \@yarg \fi
\ifnum \@yyarg >0 \@tempcnta\@yyarg \else \@tempcnta - \@yyarg \fi
\ifnum\@tempcnta>6 \@badlinearg\@tempcnta0 \fi
\setbox\@linechar\hbox{\@linefnt\@getlinechar(\@xarg,\@yyarg)}%
\ifnum \@yarg >0 \let\@upordown\raise \@clnht\z@
\else\let\@upordown\lower \@clnht \ht\@linechar\fi
\@clnwd=\wd\@linechar
\if@negarg \hskip -\wd\@linechar \def\@tempa{\hskip -2\wd \@linechar}
\else \let\@tempa\relax \fi
\@whiledim \@clnwd <\@linelen \do {\@upordown\@clnht\copy\@linechar
\@tempa \advance\@clnht \ht\@linechar \advance\@clnwd \wd\@linechar}%
\advance\@clnht -\ht\@linechar \advance\@clnwd -\wd\@linechar
\@tempdima\@linelen\advance\@tempdima -\@clnwd
\@tempdimb\@tempdima\advance\@tempdimb -\wd\@linechar
\if@negarg \hskip -\@tempdimb \else \hskip \@tempdimb \fi
\multiply\@tempdima \@m\@tempcnta \@tempdima \@tempdima \wd\@linechar
\divide\@tempcnta \@tempdima \@tempdima \ht\@linechar
\multiply\@tempdima \@tempcnta \divide\@tempdima \@m \advance\@clnht
\@tempdima
\ifdim \@linelen <\wd\@linechar \hskip \wd\@linechar
\else\@upordown\@clnht\copy\@linechar\fi}
\def\@hline{\ifnum \@xarg <0 \hskip -\@linelen \fi
\vrule height \@halfwidth depth \@halfwidth width \@linelen
\ifnum \@xarg <0 \hskip -\@linelen \fi}
\def\@getlinechar(#1,#2){\@tempcnta#1\relax
\multiply\@tempcnta 8\advance\@tempcnta -9
\ifnum #2>0 \advance\@tempcnta #2\relax
  \else\advance\@tempcnta -#2\relax\advance\@tempcnta 64 \fi
\char\@tempcnta}
\def\vector(#1,#2)#3{\@xarg #1\relax \@yarg #2\relax
\@linelen=#3\unitlength
\ifnum\@xarg =0 \@vvector \else \ifnum\@yarg =0 \@hvector \else
\@svector\fi \fi}
\def\@hvector{\@hline\hbox to 0pt{\@linefnt \ifnum \@xarg <0
\@getlarrow(1,0)\hss\else \hss\@getrarrow(1,0)\fi}}
\def\@vvector{\ifnum \@yarg <0 \@downvector \else \@upvector \fi}
\def\@svector{\@sline\@tempcnta\@yarg \ifnum\@tempcnta <0
\@tempcnta=-\@tempcnta\fi \ifnum\@tempcnta <5 \hskip -\wd\@linechar
\@upordown\@clnht \hbox{\@linefnt \if@negarg
\@getlarrow(\@xarg,\@yyarg) \else \@getrarrow(\@xarg,\@yyarg)
\fi}\else\@badlinearg\fi}
\def\@getlarrow(#1,#2){\ifnum #2 =\z@ \@tempcnta='33\else
\@tempcnta=#1\relax\multiply\@tempcnta \sixt@@n \advance\@tempcnta -9
\@tempcntb=#2\relax \multiply\@tempcntb \tw@ \ifnum \@tempcntb >0
\advance\@tempcnta \@tempcntb\relax \else\advance\@tempcnta
-\@tempcntb\advance\@tempcnta 64 \fi\fi \char\@tempcnta}
\def\@getrarrow(#1,#2){\@tempcntb=#2\relax \ifnum\@tempcntb < 0
\@tempcntb=-\@tempcntb\relax\fi \ifcase \@tempcntb\relax
\@tempcnta='55 \or \ifnum #1<3 \@tempcnta=#1\relax\multiply\@tempcnta
24 \advance\@tempcnta -6 \else \ifnum #1=3 \@tempcnta=49
\else\@tempcnta=58 \fi\fi\or \ifnum #1<3
\@tempcnta=#1\relax\multiply\@tempcnta 24 \advance\@tempcnta -3 \else
\@tempcnta=51\fi\or \@tempcnta=#1\relax\multiply\@tempcnta \sixt@@n
\advance\@tempcnta -\tw@ \else \@tempcnta=#1\relax\multiply\@tempcnta
\sixt@@n \advance\@tempcnta 7 \fi \ifnum #2<0 \advance\@tempcnta 64
\fi \char\@tempcnta}
\def\@vline{\ifnum \@yarg <0 \@downline \else \@upline\fi}
\def\@upline{\hbox to \z@{\hskip -\@halfwidth \vrule width
\@wholewidth height \@linelen depth \z@\hss}}
\def\@downline{\hbox to \z@{\hskip -\@halfwidth \vrule width
\@wholewidth height \z@ depth \@linelen \hss}}
\def\@upvector{\@upline\setbox\@tempboxa
\hbox{\@linefnt\char'66}\raise \@linelen \hbox to\z@{\lower
\ht\@tempboxa \box\@tempboxa\hss}}
\def\@downvector{\@downline\lower \@linelen \hbox to
\z@{\@linefnt\char'77\hss}}

\thinlines

\newcount\@xarg
\newcount\@yarg
\newcount\@yyarg
\newcount\@multicnt
\newdimen\@xdim
\newdimen\@ydim
\newbox\@linechar
\newdimen\@linelen
\newdimen\@clnwd
\newdimen\@clnht
\newdimen\@dashdim
\newbox\@dashbox
\newcount\@dashcnt
\catcode`@=12

\newbox\tbox
\newbox\tboxa

\def\leftzer#1{\setbox\tbox=\hbox to 0pt{#1\hss}%
\ht\tbox=0pt \dp\tbox=0pt \box\tbox}
\def\rightzer#1{\setbox\tbox=\hbox to 0pt{\hss #1}%
\ht\tbox=0pt \dp\tbox=0pt \box\tbox}
\def\centerzer#1{\setbox\tbox=\hbox to 0pt{\hss #1\hss}%
\ht\tbox=0pt \dp\tbox=0pt \box\tbox}

\def\leftput(#1,#2)#3{\setbox\tboxa=\hbox{%
\kern #1\unitlength \raise #2\unitlength\hbox{\leftzer{#3}}}%
\ht\tboxa=0pt \wd\tboxa=0pt \dp\tboxa=0pt\box\tboxa}

\def\rightput(#1,#2)#3{\setbox\tboxa=\hbox{%
\kern #1\unitlength \raise #2\unitlength\hbox{\rightzer{#3}}}%
\ht\tboxa=0pt \wd\tboxa=0pt \dp\tboxa=0pt\box\tboxa}

\def\centerput(#1,#2)#3{\setbox\tboxa=\hbox{%
\kern #1\unitlength \raise #2\unitlength\hbox{\centerzer{#3}}}%
\ht\tboxa=0pt \wd\tboxa=0pt \dp\tboxa=0pt\box\tboxa}

\unitlength=1mm

\expandafter\ifx\csname amssym.def\endcsname\relax \else
\endinput\fi
%
\expandafter\edef\csname amssym.def\endcsname{%
        \catcode`\noexpand\@=\the\catcode`\@\space}
\catcode`\@=11
%

\def\undefine#1{\let#1\undefined}
\def\newsymbol#1#2#3#4#5{\let\next@\relax
  \ifnum#2=\@ne\let\next@\msafam@\else
  \ifnum#2=\tw@\let\next@\msbfam@\fi\fi
  \mathchardef#1="#3\next@#4#5}
\def\mathhexbox@#1#2#3{\relax
  \ifmmode\mathpalette{}{\m@th\mathchar"#1#2#3}%
  \else\leavevmode\hbox{$\m@th\mathchar"#1#2#3$}\fi}
\def\hexnumber@#1{\ifcase#1 0\or 1\or 2\or 3
\or 4\or 5\or 6\or 7\or 8\or
  9\or A\or B\or C\or D\or E\or F\fi}

\font\tenmsa=msam10
\font\sevenmsa=msam7
\font\fivemsa=msam5
\newfam\msafam
\textfont\msafam=\tenmsa
\scriptfont\msafam=\sevenmsa
\scriptscriptfont\msafam=\fivemsa
\edef\msafam@{\hexnumber@\msafam}
\mathchardef\dabar@"0\msafam@39
\def\dashrightarrow{\mathrel{\dabar@\dabar@\mathchar"0
\msafam@4B}}
\def\dashleftarrow{\mathrel{\mathchar"0\msafam@4C
\dabar@\dabar@}}

\def\ulcorner{\delimiter"4\msafam@70\msafam@70 }
\def\urcorner{\delimiter"5\msafam@71\msafam@71 }
\def\llcorner{\delimiter"4\msafam@78\msafam@78 }
\def\lrcorner{\delimiter"5\msafam@79\msafam@79 }
\def\yen{{\mathhexbox@\msafam@55}}
\def\checkmark{{\mathhexbox@\msafam@58}}
\def\circledR{{\mathhexbox@\msafam@72}}
\def\maltese{{\mathhexbox@\msafam@7A}}

\font\tenmsb=msbm10
\font\sevenmsb=msbm7
\font\fivemsb=msbm5
\newfam\msbfam
\textfont\msbfam=\tenmsb
\scriptfont\msbfam=\sevenmsb
\scriptscriptfont\msbfam=\fivemsb
\edef\msbfam@{\hexnumber@\msbfam}
\def\Bbb#1{{\fam\msbfam\relax#1}}
\def\widehat#1{\setbox\z@\hbox{$\m@th#1$}%
  \ifdim\wd\z@>\tw@ em\mathaccent"0\msbfam@5B{#1}%
  \else\mathaccent"0362{#1}\fi}
\def\widetilde#1{\setbox\z@\hbox{$\m@th#1$}%
  \ifdim\wd\z@>\tw@ em\mathaccent"0\msbfam@5D{#1}%
  \else\mathaccent"0365{#1}\fi}
\font\teneufm=eufm10
\font\seveneufm=eufm7
\font\fiveeufm=eufm5
\newfam\eufmfam
\textfont\eufmfam=\teneufm
\scriptfont\eufmfam=\seveneufm
\scriptscriptfont\eufmfam=\fiveeufm
\def\frak#1{{\fam\eufmfam\relax#1}}

\csname amssym.def\endcsname

\expandafter\ifx\csname pre amssym.tex at\endcsname\relax \else
\endinput\fi
\expandafter\chardef\csname pre amssym.tex at\endcsname=\the
\catcode`\@
\catcode`\@=11
\begingroup\ifx\undefined\newsymbol \else\def\input#1
{\endgroup}\fi
\input amssym.def \relax
\newsymbol\boxdot 1200
\newsymbol\boxplus 1201
\newsymbol\boxtimes 1202
\newsymbol\square 1003
\newsymbol\blacksquare 1004
\newsymbol\centerdot 1205
\newsymbol\lozenge 1006
\newsymbol\blacklozenge 1007
\newsymbol\circlearrowright 1308
\newsymbol\circlearrowleft 1309
\undefine\rightleftharpoons
\newsymbol\rightleftharpoons 130A
\newsymbol\leftrightharpoons 130B
\newsymbol\boxminus 120C
\newsymbol\Vdash 130D
\newsymbol\Vvdash 130E
\newsymbol\vDash 130F
\newsymbol\twoheadrightarrow 1310
\newsymbol\twoheadleftarrow 1311
\newsymbol\leftleftarrows 1312
\newsymbol\rightrightarrows 1313
\newsymbol\upuparrows 1314
\newsymbol\downdownarrows 1315
\newsymbol\upharpoonright 1316
  
\newsymbol\downharpoonright 1317
\newsymbol\upharpoonleft 1318
\newsymbol\downharpoonleft 1319
\newsymbol\rightarrowtail 131A
\newsymbol\leftarrowtail 131B
\newsymbol\leftrightarrows 131C
\newsymbol\rightleftarrows 131D
\newsymbol\Lsh 131E
\newsymbol\Rsh 131F
\newsymbol\rightsquigarrow 1320
\newsymbol\leftrightsquigarrow 1321
\newsymbol\looparrowleft 1322
\newsymbol\looparrowright 1323
\newsymbol\circeq 1324
\newsymbol\succsim 1325
\newsymbol\gtrsim 1326
\newsymbol\gtrapprox 1327
\newsymbol\multimap 1328
\newsymbol\therefore 1329
\newsymbol\because 132A
\newsymbol\doteqdot 132B
  
\newsymbol\triangleq 132C
\newsymbol\precsim 132D
\newsymbol\lesssim 132E
\newsymbol\lessapprox 132F
\newsymbol\eqslantless 1330
\newsymbol\eqslantgtr 1331
\newsymbol\curlyeqprec 1332
\newsymbol\curlyeqsucc 1333
\newsymbol\preccurlyeq 1334
\newsymbol\leqq 1335
\newsymbol\leqslant 1336
\newsymbol\lessgtr 1337
\newsymbol\backprime 1038
\newsymbol\risingdotseq 133A
\newsymbol\fallingdotseq 133B
\newsymbol\succcurlyeq 133C
\newsymbol\geqq 133D
\newsymbol\geqslant 133E
\newsymbol\gtrless 133F
\newsymbol\sqsubset 1340
\newsymbol\sqsupset 1341
\newsymbol\vartriangleright 1342
\newsymbol\vartriangleleft 1343
\newsymbol\trianglerighteq 1344
\newsymbol\trianglelefteq 1345
\newsymbol\bigstar 1046
\newsymbol\between 1347
\newsymbol\blacktriangledown 1048
\newsymbol\blacktriangleright 1349
\newsymbol\blacktriangleleft 134A
\newsymbol\vartriangle 134D
\newsymbol\blacktriangle 104E
\newsymbol\triangledown 104F
\newsymbol\eqcirc 1350
\newsymbol\lesseqgtr 1351
\newsymbol\gtreqless 1352
\newsymbol\lesseqqgtr 1353
\newsymbol\gtreqqless 1354
\newsymbol\Rrightarrow 1356
\newsymbol\Lleftarrow 1357
\newsymbol\veebar 1259
\newsymbol\barwedge 125A
\newsymbol\doublebarwedge 125B
\undefine\angle
\newsymbol\angle 105C
\newsymbol\measuredangle 105D
\newsymbol\sphericalangle 105E
\newsymbol\varpropto 135F
\newsymbol\smallsmile 1360
\newsymbol\smallfrown 1361
\newsymbol\Subset 1362
\newsymbol\Supset 1363
\newsymbol\Cup 1264
  
\newsymbol\Cap 1265
  
\newsymbol\curlywedge 1266
\newsymbol\curlyvee 1267
\newsymbol\leftthreetimes 1268
\newsymbol\rightthreetimes 1269
\newsymbol\subseteqq 136A
\newsymbol\supseteqq 136B
\newsymbol\bumpeq 136C
\newsymbol\Bumpeq 136D
\newsymbol\lll 136E
  
\newsymbol\ggg 136F
  
\newsymbol\circledS 1073
\newsymbol\pitchfork 1374
\newsymbol\dotplus 1275
\newsymbol\backsim 1376
\newsymbol\backsimeq 1377
\newsymbol\complement 107B
\newsymbol\intercal 127C
\newsymbol\circledcirc 127D
\newsymbol\circledast 127E
\newsymbol\circleddash 127F
\newsymbol\lvertneqq 2300
\newsymbol\gvertneqq 2301
\newsymbol\nleq 2302
\newsymbol\ngeq 2303
\newsymbol\nless 2304
\newsymbol\ngtr 2305
\newsymbol\nprec 2306
\newsymbol\nsucc 2307
\newsymbol\lneqq 2308
\newsymbol\gneqq 2309
\newsymbol\nleqslant 230A
\newsymbol\ngeqslant 230B
\newsymbol\lneq 230C
\newsymbol\gneq 230D
\newsymbol\npreceq 230E
\newsymbol\nsucceq 230F
\newsymbol\precnsim 2310
\newsymbol\succnsim 2311
\newsymbol\lnsim 2312
\newsymbol\gnsim 2313
\newsymbol\nleqq 2314
\newsymbol\ngeqq 2315
\newsymbol\precneqq 2316
\newsymbol\succneqq 2317
\newsymbol\precnapprox 2318
\newsymbol\succnapprox 2319
\newsymbol\lnapprox 231A
\newsymbol\gnapprox 231B
\newsymbol\nsim 231C
\newsymbol\ncong 231D
\newsymbol\diagup 201E
\newsymbol\diagdown 201F
\newsymbol\varsubsetneq 2320
\newsymbol\varsupsetneq 2321
\newsymbol\nsubseteqq 2322
\newsymbol\nsupseteqq 2323
\newsymbol\subsetneqq 2324
\newsymbol\supsetneqq 2325
\newsymbol\varsubsetneqq 2326
\newsymbol\varsupsetneqq 2327
\newsymbol\subsetneq 2328
\newsymbol\supsetneq 2329
\newsymbol\nsubseteq 232A
\newsymbol\nsupseteq 232B
\newsymbol\nparallel 232C
\newsymbol\nmid 232D
\newsymbol\nshortmid 232E
\newsymbol\nshortparallel 232F
\newsymbol\nvdash 2330
\newsymbol\nVdash 2331
\newsymbol\nvDash 2332
\newsymbol\nVDash 2333
\newsymbol\ntrianglerighteq 2334
\newsymbol\ntrianglelefteq 2335
\newsymbol\ntriangleleft 2336
\newsymbol\ntriangleright 2337
\newsymbol\nleftarrow 2338
\newsymbol\nrightarrow 2339
\newsymbol\nLeftarrow 233A
\newsymbol\nRightarrow 233B
\newsymbol\nLeftrightarrow 233C
\newsymbol\nleftrightarrow 233D
\newsymbol\divideontimes 223E
\newsymbol\varnothing 203F
\newsymbol\nexists 2040
\newsymbol\Finv 2060
\newsymbol\Game 2061
\newsymbol\mho 2066
\newsymbol\eth 2067
\newsymbol\eqsim 2368
\newsymbol\beth 2069
\newsymbol\gimel 206A
\newsymbol\daleth 206B
\newsymbol\lessdot 236C
\newsymbol\gtrdot 236D
\newsymbol\ltimes 226E
\newsymbol\rtimes 226F
\newsymbol\shortmid 2370
\newsymbol\shortparallel 2371
\newsymbol\smallsetminus 2272
\newsymbol\thicksim 2373
\newsymbol\thickapprox 2374
\newsymbol\approxeq 2375
\newsymbol\succapprox 2376
\newsymbol\precapprox 2377
\newsymbol\curvearrowleft 2378
\newsymbol\curvearrowright 2379
\newsymbol\digamma 207A
\newsymbol\varkappa 207B
\newsymbol\Bbbk 207C
\newsymbol\hslash 207D
\undefine\hbar
\newsymbol\hbar 207E
\newsymbol\backepsilon 237F
\catcode`\@=\csname pre amssym.tex at\endcsname

\magnification=1200
\hsize=160 true mm
\vsize=240 true mm
\hoffset=-2mm
\voffset=8mm
\parindent=12pt   \parskip=2pt
\hfuzz=1pt
\pageno=1
\headline={\ifnum\pageno=1\else\tenrm\hfil 906-\number\folio\hfil\fi}
\footline={\hfil}

\noindent\vtop to 65 true mm{%
\hbox to 160 true mm{S\'eminaire Bourbaki\hfil Juin 2002}
\hbox to 160 true mm{54\`{e}me ann\'{e}e, 2001-2002,
n\raise4pt\hbox{o} 906\hss}
\vfil
\centerline{\bf TRAVAUX DE FRENKEL, GAITSGORY ET VILONEN}
\bigskip
\centerline{\bf  SUR LA CORRESPONDANCE DE DRINFELD-LANGLANDS}
\vfil
\centerline{\bf {\rm par} G\'{e}rard LAUMON}
\vfil}\par
\topskip=1.5 true cm
\baselineskip=13pt

En 1967, R.~Langlands a propos\'{e} une vaste extension de la
th\'{e}orie du corps de classes ab\'{e}lien de E.~Artin et J.~Tate.
Plus pr\'{e}cis\'{e}ment, il a conjectur\'{e} une correspondance
naturelle entre repr\'{e}sentations automorphes d'un groupe
r\'{e}ductif $G$ sur un corps global $F$ et repr\'{e}sentations
galoisiennes de $F$ \`{a} valeurs dans le groupe alg\'{e}brique
${}^{{\rm L}}G$ dual de $G$.  La composante neutre $\widehat{G}$ de
${}^{{\rm L}}G$ est le groupe r\'{e}ductif complexe dont les racines
sont les co-racines de $G$ et vice-versa.

Si $G$ est le groupe lin\'{e}aire $\mathop{\rm GL}(n)$ sur $F$,
$\widehat{G}$ n'est autre que $\mathop{\rm GL}(n,{\Bbb C})$ et la
corres\-pondance de Langlands globale a \'{e}t\'{e} d\'{e}montr\'{e}e
par V.~Drinfeld [Dr~1] ($n=2$) et L.~Lafforgue [La] ($n$ arbitraire)
lorsque $F$ un corps de fonctions, c'est-\`{a}-dire une extension
finie de ${\Bbb F}_{p}(t)$.  La correspondance de Langlands globale
sur les corps de nombres reste un des grands probl\`{e}mes ouverts en
math\'{e}matiques.

La {\it correspondance de Drinfeld-Langlands}, dite aussi {\it de
Langlands g\'{e}om\'{e}trique}, est un analogue conjectural de la
correspondance de Langlands pour un groupe r\'{e}ductif
d\'{e}ploy\'{e} $G$ sur une extension finie $F$ de $k(t)$, o\`{u} $k$
est un corps arbitraire.  Si $X$ est une courbe alg\'{e}brique
quasi-projective et lisse sur $k$, de corps des fonctions $F$, cette
correspondance met en dualit\'{e} un espace de modules de
$G$-fibr\'{e}s sur $X$ et un espace de modules de
$\widehat{G}$-syst\`{e}mes locaux sur $X$.

Pour $G=\mathop{\rm GL}(1)$ la correspondance de Drinfeld-Langlands
n'est autre que la th\'{e}orie du corps de classes g\'{e}om\'{e}trique
de M.~Rosenlicht et S.~Lang, expos\'{e}e par J.-P.~Serre dans [Se].
Le cas $G=\mathop{\rm GL}(2)$ a \'{e}t\'{e} trait\'{e} par V.~Drinfeld
dans les articles [Dr~2] et [Dr~3] qui sont \`{a} l'origine de la
th\'{e}orie.

De nombreux travaux ont \'{e}t\'{e} consacr\'{e}s \`{a} divers aspects
de la correspondance de Drinfeld-Langlands, en particulier ceux de
A.~Beilinson et V.~Drinfeld ([B-D]), de A.~Braverman et D.~Gaitsgory
([B-G]), et de S.~Lysenko ([Ly~1] et [Ly~2]).  Dans cet expos\'{e} je
n'\'{e}voquerai que les travaux r\'{e}cents de E.~Frenkel,
D.~Gaitsgory et K.~Vilonen dans le cas o\`{u} $X$ est projective
(correspondance partout non ramifi\'{e}e) et $G=\mathop{\rm GL}(n)$,
travaux qui g\'{e}n\'{e}ralisent ceux de V.~Drinfeld dans [Dr~2].
\vskip 2mm

Je remercie S.~Lysenko pour son aide dans la pr\'{e}paration de cet
expos\'{e}.
\vskip 10mm

\centerline{\bf 0. PR\'{E}LIMINAIRES}
\vskip 5mm

Dans le cas partout non ramifi\'{e} qui fait l'objet de cet
expos\'{e}, on esp\`{e}re \'{e}tablir une correspondance de
Drinfeld-Langlands pour chaque triplet form\'{e} d'un corps de base
$k$, d'un corps des coefficients $C$ et d'une th\'{e}orie
cohomologique \`{a} coefficients dans $C$ pour la cat\'{e}gorie des
sch\'{e}mas de type fini sur $k$.

Les trois triplets principaux (corps de base, corps des coefficients,
th\'{e}orie cohomologique) sont:
\parindent=30pt
\vskip 1mm

\itemitem{(Betti)} $k={\Bbb C}$, $C$ alg\'{e}briquement clos de
caract\'{e}ristique $0$ et la th\'{e}orie des faisceaux constructibles
de $C$-espaces vectoriels pour la topologie classique,
\vskip 1mm

\itemitem{(De Rham)} $k$ de caract\'{e}ristique nulle, $C=k$ et la
th\'{e}orie des ${\cal D}$-Modules holonomes pour la topologie de
Zariski,
\vskip 1mm

\itemitem{($\ell$-adique)} $k$ arbitraire, $C=\overline{{\Bbb
Q}}_{\ell}$ pour un nombre premier $\ell$ inversible dans $k$ et la
th\'{e}orie des faisceaux $\ell$-adiques pour la topologie \'{e}tale.
\vskip 1mm
\parindent=12pt

La donn\'{e}e premi\`{e}re de la correspondance de Drinfeld-Langlands
est celle d'une courbe alg\'{e}brique sur le corps de base $k$.  \`{A}
cette courbe on attache des espaces de modules qui sont en
g\'{e}n\'{e}ral des champs alg\'{e}briques sur $k$.  J'utiliserai
donc librement le langage des champs alg\'{e}briques (cf. [L-M]).

Les th\'{e}ories cohomologiques ci-dessus n'ont pas \'{e}t\'{e}
d\'{e}velopp\'{e}es de mani\`{e}re syst\'{e}\-ma\-ti\-que pour la
cat\'{e}gorie des champs alg\'{e}briques et les r\'{e}f\'{e}rences
sont parcellaires.  La situation est assez satisfaisante pour le
triplet ($\ell$-adique) utilis\'{e} par Frenkel, Gaitsgory et Vilonen
dans leurs articles (en fait, ils supposent de plus que $k$ est de
caract\'{e}ristique $p>0$ pour disposer d'une transformation de
Fourier g\'{e}om\'{e}trique, mais cette restriction n'est pas
n\'{e}cessaire).  Le triplet (De Rham) est utilis\'{e} dans les travaux
[B-D] de Beilinson et Drinfeld o\`{u} on trouvera une d\'{e}finition
de l'anneau des op\'{e}rateurs diff\'{e}rentiels pour un champ
alg\'{e}brique.  Le triplet (Betti), qui est en principe le plus
\'{e}l\'{e}mentaire, a \'{e}t\'{e} \'{e}tudi\'{e} par Bernstein et
Luntz [B-L], mais la th\'{e}orie des champs analytiques reste \`{a}
\'{e}crire.

Dans cet expos\'{e}, j'utiliserai n\'{e}anmoins ce dernier triplet
(Betti).  Le corps de base sera donc $k={\Bbb C}$ et, pour tout champ
alg\'{e}brique ${\cal S}$ de type fini, ou plus g\'{e}n\'{e}ralement
localement de type fini (sur ${\Bbb C}$), je noterai simplement
$D_{{\rm c}}^{{\rm b}}({\cal S})$ la cat\'{e}gorie d\'{e}riv\'{e}e des
complexes de faisceaux de $C$-espaces vectoriels \`{a} cohomologie
born\'{e}e et constructible sur le champ analytique associ\'{e} \`{a}
${\cal S}$.  Lorsque ${\cal S}$ est pr\'{e}sent\'{e} comme un quotient
d'un sch\'{e}ma $S$ par l'action d'un groupe alg\'{e}brique $G$, ce
champ analytique n'est autre que le champ quotient de $S({\Bbb
C})^{{\rm an}}$ par l'action de $G({\Bbb C})^{{\rm an}}$ et $D_{{\rm
c}}^{{\rm b}}({\cal S})$ est la cat\'{e}gorie d\'{e}riv\'{e}e
$D_{G({\Bbb C})^{{\rm an}}}^{{\rm b}} (S({\Bbb C})^{{\rm an}})$
introduite par Bernstein et Luntz dans [B-L].

Pour tous les champs alg\'{e}briques ${\cal S}$ consid\'{e}r\'{e}s
dans ce texte, la cat\'{e}gorie d\'{e}riv\'{e}e $D_{{\rm c}}^{{\rm
b}}({\cal S})$ est munie d'un produit tensoriel et d'un foncteur de
dualit\'{e} $D:D_{{\rm c}}^{{\rm b}}({\cal S})^{{\rm opp}}\rightarrow
D_{{\rm c}}^{{\rm b}}({\cal S})$; pour tous les morphismes
repr\'{e}sentables $f:{\cal S}\rightarrow {\cal T}$ consid\'{e}r\'{e}s
entre tels champs, on a des foncteurs images directes
$f_{\ast},f_{!}:D_{{\rm c}}^{{\rm b}}({\cal S}) \rightarrow D_{{\rm
c}}^{{\rm b}}({\cal T})$ et images inverses $f^{\ast},f^{!}:D_{{\rm
c}}^{{\rm b}}({\cal T})\rightarrow D_{{\rm c}}^{{\rm b}}({\cal S})$
satisfaisant au formalisme des six op\'{e}rations de Grothendieck.  La
cat\'{e}gorie $D_{{\rm c}}^{{\rm b}}({\cal S})$ est munie de la
$t$-structure pour la perversit\'{e} interm\'{e}diaire dont le
c{\oe}ur est la cat\'{e}gorie $\mathop{\rm Perv}({\cal S})$ des
faisceaux pervers sur ${\cal S}$ (cf.  [B-B-D]).  Rappelons que cette
sous-cat\'{e}gorie pleine de $D_{{\rm c}}^{{\rm b}}({\cal S})$ est
ab\'{e}lienne, noeth\'{e}rienne, artinienne et auto-duale pour la
dualit\'{e} $D$.

Si $\pi : {\cal V}\rightarrow {\cal S}$ est un fibr\'{e} vectoriel de
rang constant $r$, on peut consid\'{e}rer le ${\cal S}$-champ quotient
$\overline{{\cal V}}$ de ${\cal V}$ par l'action par homoth\'{e}tie du
groupe multiplicatif ${\Bbb G}_{{\rm m}}$.  Si ${\cal V}^{\circ}$ est
l'ouvert compl\'{e}mentaire dans ${\cal V}$ de la section nulle,
l'inclusion ${\cal V}^{\circ}\subset {\cal V}$ induit une immersion
ouverte
$$
j:{\Bbb P}({\cal V})=[{\cal V}^{\circ}/{\Bbb G}_{{\rm
m}}]=\overline{{\cal V}}^{\,\circ} \hookrightarrow \overline{{\cal V}}
$$
o\`{u} ${\Bbb P}({\cal V})\rightarrow {\cal S}$ est le fibr\'{e}
projectif des droites de ${\cal V}\rightarrow {\cal S}$.  Le ferm\'{e}
compl\'{e}mentaire de cette immersion ouverte est le quotient de la
section nulle de ${\cal V}$ par l'action triviale de ${\Bbb G}_{{\rm
m}}$, c'est-\`{a}-dire le champ classifiant $B({\Bbb G}_{{\rm
m}}/{\cal S})$.

Les complexes ${\Bbb G}_{{\rm m}}$-\'{e}quivariants constructibles de
$C$-espaces vectoriels sur ${\cal V}$ sont par d\'{e}finition les
objets de $D_{{\rm c}}^{{\rm b}}(\overline{{\cal V}})$.

On a une transformation de Fourier g\'{e}om\'{e}trique pour ces
complexes ${\Bbb G}_{{\rm m}}$-\'{e}quivariants, appel\'{e}e la {\it
transformation de Fourier homog\`{e}ne} associ\'{e}e au fibr\'{e}
vectoriel ${\cal V}/{\cal S}$ (cf.  [Lau~3]).  Cette transformation de
Fourier est une \'{e}quivalence de cat\'{e}gories d\'{e}riv\'{e}es
$$
\mathop{\rm Four}\nolimits_{\overline{{\cal V}}/{\cal S}}: D_{{\rm
c}}^{{\rm b}}(\overline{{\cal V}})\rightarrow D_{{\rm c}}^{{\rm b}}
(\overline{{\cal V}}^{\vee}),
$$
d'inverse $\mathop{\rm Four}\nolimits_{\overline{{\cal
V}}^{\vee}/{\cal S}}$, o\`{u} $\overline{{\cal V}}^{\vee}$ est le
${\cal S}$-champ quotient du fibr\'{e} vectoriel dual
$\pi^{\vee}:{\cal V}^{\vee}\rightarrow {\cal S}$ par l'action par
homoth\'{e}tie du groupe multiplicatif ${\Bbb G}_{{\rm m}}$.  Elle
commute \`{a} la dualit\'{e} et elle est $t$-exacte.  On peut la
d\'{e}finir de la fa\c{c}on suivante.

Sur le champ quotient ${\cal A}=[{\Bbb A}^{1}/{\Bbb G}_{{\rm m}}]$
quotient de la droite affine par l'action par homoth\'{e}tie du groupe
multiplicatif on a le complexe $\Psi=\beta_{\ast}C\in
\mathop{\rm ob}D_{{\rm c}}^{{\rm b}}({\cal A})$, o\`{u} $\beta
:\mathop{\rm Spec}({\Bbb C})=[{\Bbb G}_{{\rm m}}/{\Bbb G}_{{\rm
m}}]\hookrightarrow {\cal A}$ est l'immersion ouverte induite par
l'inclusion ${\Bbb G}_{{\rm m}}\subset {\Bbb A}^{1}$.  Ce complexe
admet pour faisceaux de cohomologie non triviaux
$$
{\cal H}^{0}(\Psi)=C
$$
et
$$
{\cal H}^{1}(\Psi )=\alpha_{\ast}C
$$
o\`{u} $\alpha :B({\Bbb G}_{{\rm m}})\hookrightarrow {\cal A}$ est
l'immersion ferm\'{e}e compl\'{e}mentaire de $\beta$ induite par
l'inclusion de l'origine dans ${\Bbb A}^{1}$, et c'est en fait un
faisceau pervers non irr\'{e}ductible, extension du faisceau pervers
ponctuel $\alpha_{\ast}C[-1]$ par le faisceau pervers constant $C$.

Le morphisme d'accouplement naturel entre le fibr\'{e} vectoriel
${\cal V}$ et son fibr\'{e} dual passe au quotient en un morphisme de
champs alg\'{e}briques
$$
\mu :\overline{{\cal V}}^{\vee}\times_{{\cal S}}
\overline{{\cal V}}\rightarrow {\cal A}
$$
et, si on note $p:\overline{{\cal V}}^{\vee}\times_{{\cal
S}}\overline{{\cal V}}\rightarrow \overline{{\cal V}}$ et
$p^{\vee}:\overline{{\cal V}}^{\vee}\times_{{\cal S}}
\overline{{\cal V}}\rightarrow \overline{{\cal V}}^{\vee}$ les deux
projections canoniques, on a par d\'{e}finition
$$
\mathop{\rm Four}\nolimits_{\overline{{\cal V}}/{\cal S}}(K)=
(p^{\vee})_{!}(p^{\ast}K\otimes\mu^{\ast}\Psi )[r-1],~\forall K\in
D_{{\rm c}}^{{\rm b}}(\overline{{\cal V}}).\leqno{(0.1)}
$$
Les morphismes $p$ et $p^{\vee}$ ne sont pas repr\'{e}sentables et
$(p^{\vee})_{!}$ ne respecte pas $D_{{\rm c}}^{{\rm b}}$: il
envoie $D_{{\rm c}}^{{\rm b}}(\overline{{\cal V}}^{\vee}\times_{{\cal
S}} \overline{{\cal V}})$ dans $D_{{\rm c}}^{-}(\overline{{\cal
V}}^{\vee})$.  Cependant, on peut v\'{e}rifier que $\mathop{\rm
Four}\nolimits_{\overline{{\cal V}}/{\cal S}}$ respecte lui $D_{{\rm
c}}^{{\rm b}}$.

\thm D\'{E}FINITION (0.2)
\enonce
Un faisceau pervers $L$ sur ${\Bbb P}({\cal V})$ est dit {\rm propre}
si la fl\`{e}che d'oubli des supports $j_{!}L\rightarrow j_{\ast} L$
est un isomorphisme.  Si tel est le cas $j_{!}L\cong j_{\ast}L$ est en
fait un faisceau pervers sur $\overline{{\cal V}}$ qui n'est autre que
le prolongement interm\'{e}diaire $j_{!\ast}L$ de $L$.
\endthm

\thm LEMME (0.3)
\enonce
Soit $L$ un faisceau pervers sur ${\Bbb P}({\cal V})$.  Alors, $L$ est
propre si et seulement si $(\overline{\pi}^{\,\circ})_{!}L=0$ o\`{u}
$\overline{\pi}^{\,\circ}:{\Bbb P}({\cal V})\rightarrow {\cal S}$ est
la projection canonique.
\hfill\hfill$\square$
\endthm

La transformation de Radon g\'{e}om\'{e}trique associ\'{e}e au
fibr\'{e} projectif ${\Bbb P}({\cal V})\rightarrow {\cal S}$ est le foncteur
$$
\mathop{\rm Rad}\nolimits_{{\Bbb P}({\cal V})/{\cal S}}:D_{{\rm
c}}^{{\rm b}} ({\Bbb P}({\cal V}))\rightarrow D_{{\rm c}}^{{\rm b}}
({\Bbb P}({\cal V}^{\vee}))
$$
d\'{e}fini par
$$
\mathop{\rm Rad}\nolimits_{{\Bbb P}({\cal V})/{\cal S}}(L)=
(q^{\vee})_{!}q^{\ast}L[r-1]
$$
o\`{u} $q:{\cal H}\rightarrow {\Bbb P}({\cal V})$ et $q^{\vee}:{\cal
H}\rightarrow {\Bbb P}({\cal V}^{\vee})$ sont les deux projections
canoniques du champ d'incidence
$$
{\cal H}\subset {\Bbb P}({\cal V}^{\vee})\times_{{\cal S}}{\Bbb P}
({\cal V})
$$
des couples form\'{e}s d'une droite et d'un hyperplan de ${\cal V}
\rightarrow {\cal S}$ tels que la droite soit contenue dans
l'hyperplan.  On renvoie \`{a} la monographie de Brylinski [Br] pour
une \'{e}tude d\'{e}taill\'{e}e de cette transformation.

\thm LEMME (0.4)
\enonce
Soit $L$ un faisceau pervers sur ${\Bbb P}({\cal V})$ que l'on suppose
propre.  Alors $\mathop{\rm Four} \nolimits_{\overline{{\cal V}}/{\cal
S}}(j_{!}L)$ est aussi propre et sa restriction \`{a} l'ouvert ${\Bbb
P}({\cal V}^{\vee})= \overline{{\cal V}}^{\vee\circ}$
compl\'{e}mentaire de la section nulle dans $\overline{{\cal
V}}^{\vee}$ est \'{e}gale \`{a} $\mathop{\rm Rad}\nolimits_{{\Bbb
P}({\cal V})/{\cal S}}(L)$.
\hfill\hfill$\square$
\endthm
\vskip 5mm

\centerline{\bf 1. CORPS DE CLASSES G\'{E}OM\'{E}TRIQUE}
\vskip 5mm

Soit $X$ une courbe alg\'{e}brique connexe, projective et lisse (une
surface de Riemann compacte connexe) de genre $g$.  Fixons un point
base $x_{0}$, une base $(\delta_{1},\ldots ,\delta_{2g})$ du ${\Bbb
Z}$-module libre
$$
H_{1}(X,{\Bbb Z})=\pi_{1}^{{\rm ab}}(X):=\pi_{1}(X,x_{0})/[\pi_{1}(X,x_{0}),
\pi_{1}(X,x_{0})]
$$
et une base $(\omega_{1},\ldots ,\omega_{g})$ de l'espace vectoriel
complexe $H^{0}(X,\Omega_{X}^{1})$.  Les p\'{e}riodes
$$
\Pi_{i}=\left(\int_{\delta_{i}}\omega_{1},\ldots
,\int_{\delta_{i}}\omega_{g}\right)\in {\Bbb C}^{g},~i=1,\ldots ,2g,
$$
sont lin\'{e}airement ind\'{e}pendantes sur ${\Bbb Z}$ et engendrent
un r\'{e}seau
$\Lambda\subset {\Bbb C}^{g}$.  La jacobienne de $X$ est le tore
complexe de dimension $g$,
$$
J(X)={\Bbb C}^{g}/\Lambda .
$$
On a une application analytique
$$
\varphi : X\rightarrow J(X),~x\mapsto\left(\int_{x_{0}}^{x}
\omega_{1},\ldots ,\int_{x_{0}}^{x}\omega_{g}\right)+\Lambda ,
$$
qui envoie $x_{0}$ sur l'\'{e}l\'{e}ment neutre $0$ de $J(X)$.

Soit $\mathop{\rm Pic}(X)$ le sch\'{e}ma de Picard qui param\`{e}tre
les classes d'isomorphie de fibr\'{e}s en droites sur $X$.  C'est un
sch\'{e}ma en groupes pour la structure de groupe induite par le
produit tensoriel des fibr\'{e}s en droites, qui s'ins\`{e}re dans la
suite exacte
$$
1\rightarrow \mathop{\rm Pic}\nolimits^{0}(X)\rightarrow \mathop{\rm
Pic}(X)\,\smash{\mathop{\hbox to 8mm{\rightarrowfill}}
\limits^{\scriptstyle \mathop{\rm deg}}}\,{\Bbb Z}\rightarrow 0
$$
o\`{u} $\mathop{\rm deg}({\cal L})$ est le degr\'{e} du fibr\'{e} en
droites ${\cal L}$, et sa composante neutre $\mathop{\rm Pic}^{0}(X)$
est une vari\'{e}t\'{e} ab\'{e}lienne.

\thm TH\'{E}OR\`{E}ME 1.1 (Abel-Jacobi)
\enonce
La jacobienne $J(X)$ est le tore complexe sous-jacent \`{a} la
vari\'{e}t\'{e} ab\'{e}lienne $\mathop{\rm Pic}\nolimits^{0}(X)$ et le
morphisme analytique
$$
\varphi :X\rightarrow J(X)\cong\mathop{\rm Pic}\nolimits^{0}(X)
$$
n'est autre que le morphisme alg\'{e}brique qui envoie $x\in X$ sur la
classe d'isomorphie de ${\cal O}_{X}([x]-[x_{0}])$ dans $\mathop{\rm
Pic}\nolimits^{0}(X)$.
\hfill\hfill$\square$
\endthm

L'homomorphisme
$$
\pi_{1}(X,x_{0})\rightarrow\pi_{1}(J(X),0)
$$
induit par $\varphi$ est l'application quotient
$$
\pi_{1}(X,x_{0})\twoheadrightarrow
\pi_{1}(X,x_{0})/[\pi_{1}(X,x_{0}),\pi_{1}(X,x_{0})]=\Lambda .
$$
Pour tout caract\`{e}re $\chi :\pi_{1}(X,x_{0})\rightarrow C^{\times}$
il existe donc un unique caract\`{e}re $\mathop{\rm Aut}\nolimits_{\chi}:
\pi_{1}(\mathop{\rm Pic}\nolimits^{0}(X),0)\rightarrow C^{\times}$
tel que $\chi =\mathop{\rm Aut}\nolimits_{\chi}\circ\varphi_{\ast}$.

Comme un syst\`{e}me local $E$ (de $C$-espaces vectoriels) sur une
vari\'{e}t\'{e} n'est autre qu'une repr\'{e}sentation de
dimension finie sur $C$ du groupe fondamental de cette vari\'{e}t\'{e}, on a
montr\'{e} par voie analytique :

\thm TH\'{E}OR\`{E}ME 1.2
\enonce
Pour tout syst\`{e}me local $E$ de rang $1$ sur $X$, il existe un
unique syst\`{e}me local {\rm (}rigidifi\'{e} \`{a} l'origine{\rm )}
$\mathop{\rm Aut}\nolimits_{E}$ de rang $1$ sur $\mathop{\rm
Pic}\nolimits^{0}(X)$ tel que $\varphi^{\ast}\mathop{\rm
Aut}\nolimits_{E}=E$.

De plus, si $m:\mathop{\rm Pic}\nolimits^{0}(X)\times \mathop{\rm Pic}
\nolimits^{0}(X)\rightarrow \mathop{\rm Pic}\nolimits^{0}(X)$ est la
loi de groupe, on a un isomorphisme canonique
$$
m^{\ast}\mathop{\rm Aut}\nolimits_{E}\cong \mathop{\rm
Aut}\nolimits_{E}\boxtimes \mathop{\rm Aut}\nolimits_{E}.
$$
\endthm

\rem D\'{e}monstration alg\'{e}brique du th\'{e}or\`{e}me $1.2$
\endrem
Pour chaque entier $d\geq 0$, le groupe sym\'{e}trique ${\frak S}_{d}$ agit
de mani\`{e}re \'{e}vidente sur le produit
$$
X^{d}=\overbrace{X\times\cdots\times X}^{d}.
$$
Le sch\'{e}ma quotient
$$
X^{(d)}=X^{d}/{\frak S}_{d},
$$
qui est connexe, projectif et lisse de dimension $d$, est l'espace de
modules des diviseurs $D=\sum_{x}d_{x}[x]$ effectifs ($d_{x}\geq 0$)
et de degr\'{e} $\sum_{x}d_{x}=d$ sur $X$, et le morphisme $\varphi$
induit un morphisme
$$
\varphi^{(d)}:X^{(d)}\rightarrow \mathop{\rm Pic}
\nolimits^{0}(X),~D\mapsto {\cal O}_{X}(D-d[x_{0}]).
$$
D\`{e}s que $d\geq 2g-1$, $\varphi^{(d)}$ est un fibr\'{e} projectif
de rang $d-g$ et a donc ses fibres simplement connexes.

Pour chaque syst\`{e}me local $E$ de rang $n$ sur $X$ et chaque entier
$d\geq 0$, le produit tensoriel externe $E^{\boxtimes d}$ est un
syst\`{e}me local de rang $n^{d}$ sur $X^{d}$ muni d'une action de
${\frak S}_{d}$ qui rel\`{e}ve celle sur $X^{d}$.  Si on note
$r:X^{d}\twoheadrightarrow X^{d}/{\frak S}_{d}=X^{(d)}$ le morphisme
quotient, le faisceau constructible de $C$-espaces vectoriels
$r_{\ast}E^{\boxtimes d}$ est donc muni d'une action de ${\frak
S}_{d}$ et on note
$$
E^{(d)}=(r_{\ast}E^{\boxtimes d})^{{\frak S}_{d}}\subset E^{\boxtimes
d}
$$
le sous-faisceau des vecteurs fixes par cette action.  On v\'{e}rifie
que la fibre en $D\in X^{(d)}$ de $E^{(d)}$ est \'{e}gale \`{a}
$$
(E^{(d)})_{D}=\bigotimes_{x}\mathop{\rm Sym}\nolimits^{d_{x}}
E_{x}.\leqno{(1.3)}
$$

Si maintenant $E$ est de rang $1$ sur $X$, $E^{(d)}$ est en fait un
syst\`{e}me local de rang $1$ sur $X^{(d)}$.  Pour tout $d\geq 2g-1$,
$E^{(d)}$ est donc constant sur les fibres simplement connexes du
morphisme $\varphi^{(d)}:X^{(d)}\rightarrow \mathop{\rm Pic}
\nolimits^{0}(X)$ et se descend en un syst\`{e}me local sur
$\mathop{\rm Pic}\nolimits_{}^{0}(X)$ qui n'est autre que $\mathop{\rm
Aut}\nolimits_{E}$.
\hfill\hfill$\square$
\vskip 3mm

C'est le cas $\mathop{\rm GL}(1)$ de la correspondance de
Drinfeld-Langlands partout non ramifi\'{e}e.
\vskip 5mm

\centerline{\bf 2. LE TH\'{E}OR\`{E}ME PRINCIPAL}
\vskip 5mm

Soit $n$ un entier $\geq 2$.  La g\'{e}n\'{e}ralisation naturelle de
$\mathop{\rm Pic}(X)$, ou plut\^{o}t du champ quotient du sch\'{e}ma
$\mathop{\rm Pic}(X)$ par l'action triviale de ${\Bbb G}_{{\rm m}}$,
est le champ $\mathop{\rm Fib}\nolimits_{n}(X)$ des fibr\'{e}s
vectoriels de rang $n$ sur la surface de Riemann $X$.  Ce champ est
alg\'{e}brique, lisse purement de dimension $n^{2}(g-1)$ et
r\'{e}union croissante d'ouverts de type fini.  Si $g\geq 2$, il
contient un ouvert dense qui est une ${\Bbb G}_{{\rm m}}$-gerbe sur le
sch\'{e}ma quasi-projectif de modules des fibr\'{e}s stables sur $X$.
Les composantes connexes $\mathop{\rm Fib}\nolimits_{n}^{d}(X)$ de
$\mathop{\rm Fib}\nolimits_{n}(X)$ sont d\'{e}coup\'{e}es par le
degr\'{e} $d$ du fibr\'{e} vectoriel universel.

Un fibr\'{e} vectoriel ${\cal L}$ de rang $n$ sur $X$ peut aussi
\^{e}tre vu comme un ${\cal O}_{X}$-Module localement libre de rang
$n$.  Une {\it modification inf\'{e}rieure \'{e}l\'{e}mentaire} de ${\cal
L}$ en un point $x\in X$ est un fibr\'{e} vectoriel ${\cal L}'$ de
rang $n$ sur $X$ qui est contenu dans ${\cal L}$ en tant que
sous-${\cal O}_{X}$-Module de telle sorte que le ${\cal O}_{X}$-Module
coh\'{e}rent ${\cal L}/{\cal L}'$ soit un faisceau gratte-ciel de
longueur $1$ concentr\'{e} en $x$.  Le degr\'{e} de ${\cal L}'$ est
donc \'{e}gal \`{a}
$$
\mathop{\rm deg}({\cal L}')=\mathop{\rm deg}({\cal L})-1.
$$

Le {\it champ de Hecke} est le champ modulaire ${\mathop{\rm
Hecke}\nolimits_{n}(X)}$ des triplets $(x,{\cal L},{\cal L}'\subset
{\cal L})$ o\`{u} $x\in X$, ${\cal L}\in \mathop{\rm
Fib}\nolimits_{n}(X)$ et ${\cal L}'\subset {\cal L}$ est une
modification \'{e}l\'{e}mentaire inf\'{e}rieure de ${\cal L}$ en $x$.

La {\it correspondance de Hecke} est le diagramme
$$\diagram{
\noalign{\vskip 10mm}
\put (16,10){\mathop{\rm Hecke}\nolimits_{n}(X)}
\arrow(14,8)\dir(-1,-1)\length{8}
\put (4,5){\scriptstyle (q,p')}\cr
\arrow(33,8)\dir(1,-1)\length{8}\cr
\put (38,5){\scriptstyle p}\cr
\kern -5mm X\times \mathop{\rm Fib}\nolimits_{n}(X)&&
\kern 14mm\mathop{\rm Fib}\nolimits_{n}(X)\cr}
$$
o\`{u} les projections $(q,p')$ et $p$ sont donn\'{e}es par
$$
(q,p')(x,{\cal L},{\cal L}'\subset{\cal L})=(x,{\cal L}')
$$
et
$$
p(x,{\cal L},{\cal L}'\subset{\cal L})={\cal L},
$$
et sont des morphismes repr\'{e}sentables projectifs et lisses de
dimension relative $n-1$ et $n$ respectivement.

Cette correspondance induit un  op\'{e}rateur, dit {\it de Hecke}, sur la
cat\'{e}gorie d\'{e}riv\'{e}e $D_{{\rm c}}^{{\rm b}}(\mathop{\rm Fib}
\nolimits_{n}(X))$, qui envoie $K\in D_{{\rm c}}^{{\rm
b}}(\mathop{\rm Fib} \nolimits_{n}(X))$ sur
$$
H(K)=(q,p')_{\ast}p^{\ast}K[n-1]\in D_{{\rm c}}^{{\rm b}}
(X\times\mathop{\rm Fib}\nolimits_{n}(X)).
$$
On peut {\og}{it\'{e}rer}{\fg} $H$ et on obtient en particulier
l'op\'{e}rateur
$$
H\circ H:D_{{\rm c}}^{{\rm b}}(\mathop{\rm Fib}\nolimits_{n}(X))
\rightarrow D_{{\rm c}}^{{\rm b}}(X\times X\times\mathop{\rm
Fib}\nolimits_{n}(X)),
$$
dont la restriction au compl\'{e}mentaire $(X\times
X-\Delta_{X})\times \mathop{\rm Fib}\nolimits_{n}(X)$ de la diagonale
est ${\frak S}_{2}$-\'{e}quivariante pour l'action qui permute les
deux copies de $X$.

\thm D\'{E}FINITION 2.1
\enonce
Soit $E$ un syst\`{e}me local irr\'{e}ductible de rang $n$ sur $X$.
Un faisceau pervers irr\'{e}ductible $K$ sur $\mathop{\rm Fib}
\nolimits_{n}(X)$ est dit avoir la {\rm propri\'{e}t\'{e} de Hecke
relativement \`{a}} $E$ s'il existe un isomorphisme
$$
H(K)\cong E\boxtimes K
$$
tel que la restriction de l'isomorphisme induit
$$
(H\circ H)(K)\cong E\boxtimes E\boxtimes K
$$
\`{a} l'ouvert $(X\times X-\Delta_{X})\times \mathop{\rm Fib}
\nolimits_{n}(X)$ soit ${\frak S}_{2}$-\'{e}quivariant.
\endthm

\thm TH\'{E}OR\`{E}ME 2.2 (Drinfeld pour $n=2$, [Dr~2]; Frenkel,
Gaitsgory et Vilonen en g\'{e}n\'{e}ral, [F-G-V~2] et [Ga])
\enonce
Pour tout syst\`{e}me local irr\'{e}ductible $E$ de rang $n$ sur $X$,
il existe un faisceau pervers $\mathop{\rm Aut}\nolimits_{E}$ sur
$\mathop{\rm Fib} \nolimits_{n}(X)$ dont la restriction \`{a} chaque
composante connexe $\mathop{\rm Fib}\nolimits_{n}^{d}(X)$ est
irr\'{e}ductible et qui a la propri\'{e}t\'{e} de Hecke relativement
\`{a} $E$.
\endthm

La d\'{e}monstration du th\'{e}or\`{e}me se fait en deux temps.

Premi\`{e}rement on construit un complexe de faisceaux constructibles
$\mathop{\rm Aut}\nolimits_{E}^{\prime}$ sur l'espace de modules
$\mathop{\rm Fib}\nolimits_{n}^{\prime}(X)$ des couples $(L,s)$
form\'{e}s d'un fibr\'{e} vectoriel $L$ de rang $n$ sur $X$ et d'une
section non nulle \`{a} homoth\'{e}tie pr\`{e}s
$$
s:(\Omega_{X}^{1})^{\otimes (n-1)}\hookrightarrow L.
$$
Cette construction est inspir\'{e}e par le d\'{e}veloppement en
s\'{e}rie de Fourier d'une forme automorphe cuspidale pour
$\mathop{\rm GL}(n)$ et la formule de Shintani pour les coefficients
de cette s\'{e}rie de Fourier, compte tenu du dictionnaire
fonctions-faisceaux de Grothendieck; mais ce n'en est pas une simple
transposition.

Puis on montre que, sur un gros ouvert de $\mathop{\rm
Fib}\nolimits_{n}^{\prime}(X)$, $\mathop{\rm Aut}
\nolimits_{E}^{\prime}$ est un faisceau pervers qui se descend par le
morphisme $\mathop{\rm Fib} \nolimits_{n}^{\prime}(X)\rightarrow
\mathop{\rm Fib} \nolimits_{n}(X)$ d'oubli de la section $s$, en un
faisceau pervers $\mathop{\rm Aut}\nolimits_{E}$ sur $\mathop{\rm
Fib}\nolimits_{n}(X)$ qui a les propri\'{e}t\'{e}s requises.  C'est
bien entendu dans cette deuxi\`{e}me \'{e}tape que r\'{e}side toute la
difficult\'{e}.
\vskip 5mm

\centerline{\bf 3.  LES FAISCEAUX PERVERS ${\cal L}_{E}^{d}$}
\vskip 5mm

Pour chaque entier $d\geq 0$, soit $\mathop{\rm Coh}
\nolimits_{0}^{d}(X)$ le champ des ${\cal O}_{X}$-Modules
coh\'{e}rents ${\cal M}$ de rang g\'{e}n\'{e}rique $0$,
c'est-\`{a}-dire \`{a} support fini, et de longueur $\mathop{\rm dim}
H^{0}(X,{\cal M})=d$.  On note $\mathop{\rm Coh}\nolimits_{0}(X)$ la
r\'{e}union disjointe des $\mathop{\rm Coh}\nolimits_{0}^{d}(X)$.

La d\'{e}finition du champ $\mathop{\rm Coh} \nolimits_{0}^{d}(X)$ a
un sens pour toute courbe alg\'{e}brique $X$ quasi-projective; en
particulier, on peut consid\'{e}rer le champ $\mathop{\rm
Coh}\nolimits_{0}^{d}({\Bbb A}^{1})$.  Comme la donn\'{e}e d'un ${\cal
O}_{{\Bbb A}^{1}}$-Module coh\'{e}rent de rang g\'{e}n\'{e}rique $0$
et de longueur $d$ \'{e}quivaut \`{a} la donn\'{e}e d'un espace
vectoriel de dimension $d$ muni d'un endomorphisme, $\mathop{\rm
Coh}\nolimits_{0}^{d}({\Bbb A}^{1})$ n'est autre que le champ quotient
$[\mathop{\rm gl}(d)/\mathop{\rm GL}(d)]$ de l'espace affine des
matrices carr\'{e}es de taille $d\times d$ par l'action par
conjugaison de $\mathop{\rm GL}(d)$.  Comme toute courbe
alg\'{e}brique lisse $X$ est localement pour la topologie \'{e}tale
isomorphe \`{a} la droite affine ${\Bbb A}^{1}$, $\mathop{\rm Coh}
\nolimits_{0}^{d}(X)$ est localement pour la topologie \'{e}tale
isomorphe \`{a} $[\mathop{\rm gl}(d)/\mathop{\rm GL}(d)]$.

Chaque diviseur effectif $D$ sur $X$ d\'{e}finit un ${\cal
O}_{X}$-Module de torsion
$$
{\cal O}_{X,D}={\cal O}_{X}/{\cal O}_{X}(-D)
$$
de longueur le degr\'{e} de $D$ et, pour chaque entier $d\geq 0$, le
morphisme
$$
X^{(d)}\rightarrow \mathop{\rm Coh}\nolimits_{0}^{d}(X),~ D\mapsto
{\cal O}_{X,D},
$$
est lisse de dimension relative $d$.  Tout ${\cal M}\in\mathop{\rm
Coh}\nolimits_{0}(X)$ est isomorphe \`{a} $\bigoplus_{i\in I}{\cal
O}_{X,D_{i}}$ pour une famille finie $(D_{i})_{i\in I}$ de diviseurs
effectifs sur $X$ et, pour chaque entier $d\geq 0$, on a un morphisme
{\it d\'{e}terminant}
$$
\mathop{\rm det}:\mathop{\rm Coh}\nolimits_{0}^{d}(X)\rightarrow X^{(d)}
$$
qui envoie ${\cal M}\cong \bigoplus_{i\in I}{\cal O}_{X,D_{i}}$ sur
$\sum_{i\in I}D_{i}$.  On note $X_{{\rm rss}}^{(d)}$ l'ouvert de
$X^{(d)}$ form\'{e} des $D$ sans multiplicit\'{e} ($d_{x}=0\hbox{ ou
}1$ quel que soit $x$) et $\mathop{\rm Coh}\nolimits_{0,{\rm
rss}}^{d}(X)=\mathop{\rm det}\nolimits^{-1}(X_{{\rm rss}}^{(d)})$
l'ouvert correspondant de $\mathop{\rm Coh}\nolimits_{0}^{d}(X)$.

Soit $E$ un syst\`{e}me local de rang $n$ sur $X$.  La formule (1.3)
pour la fibre en $D\in X^{(d)}$ de $E^{(d)}$ montre que le rang de
cette fibre varie avec $D$, mais que, au-dessus de l'ouvert $X_{{\rm
rss}}^{(d)}$, ce rang est constant et la restriction de $E^{(d)}$ est
un syst\`{e}me local de rang $n^{d}$.

\thm LEMME 3.1
\enonce
Pour chaque entier $d\geq 0$, le complexe $E^{(d)}[d]\in D_{{\rm
c}}^{{\rm b}}(X^{(d)})$, constitu\'{e} du faisceau constructible
$E^{(d)}$ plac\'{e} en degr\'{e} $-d$, est un faisceau pervers, qui est
le prolongement interm\'{e}diaire de sa restriction \`{a} $X_{{\rm
rss}}^{(d)}$.

Le faisceau pervers $E^{(d)}[d]$ est irr\'{e}ductible {\rm (}resp.
semi-simple{\rm )} si $E$ l'est.
\hfill\hfill$\square$
\endthm

La restriction
$$
{\cal L}_{E,{\rm rss}}^{d}=\mathop{\rm det}\nolimits^{\ast}E^{(d)}|
\mathop{\rm Coh}\nolimits_{0, {\rm rss}}^{d}(X)
$$
du faisceau constructible $\mathop{\rm det}\nolimits^{\ast}E^{(d)}$
\`{a} l'ouvert dense $\mathop{\rm Coh}\nolimits_{0,{\rm rss}}^{d}(X)$
de $\mathop{\rm Coh}\nolimits_{0}^{d}(X)$ est aussi un syst\`{e}me
local de rang $n^{d}$, qui est irr\'{e}ductible (resp.  semi-simple)
si $E$ l'est.

\thm D\'{E}FINITION 3.2
\enonce
Pour chaque entier $d\geq 0$, le faisceau pervers ${\cal L}_{E}^{d}$
associ\'{e} \`{a} un syst\`{e}me local $E$ sur $X$ est le prolongement
interm\'{e}diaire de ${\cal L}_{E,{\rm rss}}^{d}$ \`{a} $\mathop{\rm
Coh}\nolimits_{0}^{d}(X)$ tout entier.
\endthm

Bien s\^{u}r, ${\cal L}_{E}^{d}$ est irr\'{e}ductible (resp.
semi-simple) si $E$ l'est.  On v\'{e}rifie en outre que l'image
r\'{e}ciproque d\'{e}cal\'{e}e de $d$ de ${\cal L}_{E}^{d}$ par le
morphisme $X^{(d)}\rightarrow \mathop{\rm Coh}
\nolimits_{0}^{d},~D\mapsto {\cal O}_{X,D}$, n'est autre que le
faisceau pervers $E^{(d)}[d]$.
\vskip 5mm

\centerline{\bf 4.  LE COMPLEXE $\mathop{\rm Aut}
\nolimits_{E}^{\prime}$}
\vskip 5mm

{\it On suppose dor\'{e}navant que le genre $g$ de $X$ est $\geq 2$.}
En genre $0$ ou $1$, la corres\-pondance de Drinfeld-Langlands partout
non ramifi\'{e}e pour $\mathop{\rm GL}(n)$ avec $n\geq 2$ est vide
puisqu'il n'y a pas de syst\`{e}me local irr\'{e}ductible de rang
$\geq 2$ sur $X$.
\vskip 2mm

On se propose dans cette section de construire un candidat
$\mathop{\rm Aut}\nolimits_{E}^{\prime}$ pour la restriction de
$\mathop{\rm Aut}\nolimits_{E}$ au champ alg\'{e}brique $\mathop{\rm
Fib}\nolimits_{n}^{\prime}(X)$ des fibr\'{e}s vectoriels ${\cal L}$ de
rang $n$ sur $X$ munis d'une section $s:(\Omega_{X}^{1})^{\otimes n-1}
\rightarrow {\cal L}$ \`{a} homoth\'{e}tie pr\`{e}s.
\vskip 2mm

Consid\'{e}rons le champ alg\'{e}brique ${\cal U}_{n}$ des triplets
$({\cal L},{\cal L}_{\bullet},s_{\bullet})$ o\`{u} ${\cal L}_{n}$
est un fibr\'{e} vectoriel de rang $n$,
$$
{\cal L}_{\bullet}=((0)={\cal L}_{0}\subset {\cal L}_{1}\subset
\cdots\subset {\cal L}_{n}={\cal L})
$$
est un drapeau complet de sous-fibr\'{e}s vectoriels et $s_{\bullet}$
est une famille d'isomorphismes de fibr\'{e}s en droites
$$
s_{j}:(\Omega_{X}^{1})^{\otimes j-1}\buildrel\sim\over\longrightarrow
{\cal L}_{n-j+1}/{\cal L}_{n-j},~j=1,\ldots ,n,
$$
cette famille \'{e}tant prise \`{a} homoth\'{e}tie pr\`{e}s:
$(s_{j})_{j=1,\ldots ,n}\sim (ts_{j})_{j=1,\ldots ,n}$ pour tout $t\in
{\Bbb G}_{{\rm m}}$.  On a des morphismes
$$\diagram{
{\cal U}_{n}&\kern -14mm\smash{\mathop{\hbox to 21mm{\rightarrowfill}}
\limits^{\scriptstyle g_{n}}}\kern -1mm&{\Bbb A}^{1}\cr
\llap{$\scriptstyle f_{n}$}\left\downarrow
\vbox to 4mm{}\right.\rlap{}&&\cr
\mathop{\rm Fib}\nolimits_{n}^{\prime\, n(n-1)(g-1)}(X)&&\cr}
$$
o\`{u} $f_{n}$ envoie le triplet $({\cal L},{\cal L}_{\bullet},
s_{\bullet})$ sur le fibr\'{e} ${\cal L}$ muni de la section \`{a}
homoth\'{e}tie pr\`{e}s induite par
$$
s_{n}:(\Omega_{X}^{1})^{\otimes n-1}\buildrel\sim\over\longrightarrow
{\cal L}_{1}\subset {\cal L}_{n}={\cal L}
$$
et $g_{n}$ envoie ce m\^{e}me triplet sur le scalaire somme des
classes d'extensions
$$
0\rightarrow {\cal L}_{j}/{\cal L}_{j-1}\rightarrow {\cal L}_{j+1}/{\cal
L}_{j-1}\rightarrow {\cal L}_{j+1}/{\cal L}_{j}\rightarrow 0
$$
dans
$$
\mathop{\rm Ext}\nolimits_{{\cal O}_{X}}^{1}({\cal L}_{j+1}/{\cal
L}_{j},{\cal L}_{j}/{\cal L}_{j-1})\cong \mathop{\rm Ext}
\nolimits_{{\cal O}_{X}}^{1}((\Omega_{X}^{1})^{\otimes n-j-1},
(\Omega_{X}^{1})^{\otimes n-j})\cong \mathop{\rm Ext}\nolimits_{{\cal
O}_{X}}^{1} ({\cal O}_{X},\Omega_{X}^{1})\cong {\Bbb C}.
$$
Si l'on fait agir le groupe multiplicatif ${\Bbb G}_{{\rm m}}$ sur
${\cal U}_{n}$ par
$$
t\cdot ({\cal L},{\cal L}_{\bullet},s_{\bullet})=({\cal L},{\cal
L}_{\bullet},t\cdot s_{\bullet})
$$
o\`{u} $t\cdot s_{\bullet}=(t^{n-j}s_{j})_{j=1,\ldots ,n}$ et par
homoth\'{e}tie sur ${\Bbb A}^{1}$, le morphisme $g_{n}$ est ${\Bbb
G}_{{\rm m}}$-\'{e}quivariant.  En passant au quotient par cette
action on obtient des morphismes
$$\diagram{
\overline{{\cal U}}_{n}&\kern -14mm\smash{\mathop{\hbox to
21mm{\rightarrowfill}} \limits^{\scriptstyle \overline{g}_{n}}}\kern
-1mm&[{\Bbb A}^{1}/{\Bbb G}_{{\rm m}}]={\cal A}\cr
\llap{$\scriptstyle \overline{f}_{n}$}\left\downarrow
\vbox to 4mm{}\right.\rlap{}&&\cr
\mathop{\rm Fib}\nolimits_{n}^{\prime\, n(n-1)(g-1)}(X)&&\cr}
$$
et on peut former le complexe de faisceaux constructibles
$$
L_{n}=(\overline{f}_{n})_{!}(\overline{g}_{n})^{\ast}\Psi [\mathop{\rm
dim}\overline{{\cal U}}_{n}]
$$
o\`{u} $\Psi\in D_{{\rm c}}^{{\rm b}}({\cal A})$ est le complexe
d\'{e}fini dans la section $0$.

Pour chaque entier $d$, on a une correspondance
$$\diagram{
\noalign{\vskip 10mm}
\put (25,10){\mathop{\rm Mod}\nolimits_{n}^{\prime\, d}(X)}
\put (42,10){\smash{\mathop{\hbox to 8mm{\rightarrowfill}}
\limits^{\scriptstyle q}}}
\put (52,10){\mathop{\rm Coh}\nolimits_{0}^{d}(X)}
\arrow(23,8)\dir(-1,-1)\length{8}
\put (17,5){\scriptstyle p'}\cr
\arrow(41,8)\dir(1,-1)\length{8}\cr
\put (46,5){\scriptstyle p}\cr
\kern -3mm \mathop{\rm Fib}\nolimits_{n}^{\prime\, n(n-1)(g-1)}(X)&&
\kern 3mm\mathop{\rm Fib}\nolimits_{n}^{\prime\, n(n-1)(g-1)+d}(X)\cr}
$$
o\`{u} ${\mathop{\rm Mod}\nolimits_{n}^{\prime\, d}(X)}$ est le champ
des diagrammes
$$
(\Omega_{X}^{1})^{\otimes n-1}\hookrightarrow {\cal L}'\subset {\cal
L}
$$
form\'{e}s de deux fibr\'{e}s vectoriels de rang $n$ et de degr\'{e}s
$n(n-1)(g-1)$ et $n(n-1)(g-1)+d$ respectivement et d'injections ${\cal
O}_{X}$-lin\'{e}aires $(\Omega_{X}^{1})^{\otimes n-1}\hookrightarrow
{\cal L}'$ (\`{a} homoth\'{e}tie pr\`{e}s) et ${\cal L}'
\hookrightarrow {\cal L}$, et o\`{u} $p$, $p'$ et $q$ envoient un tel
diagramme sur $(\Omega_{X}^{1})^{\otimes n-1}\hookrightarrow {\cal
L}'$ (\`{a} homoth\'{e}tie pr\`{e}s), $(\Omega_{X}^{1})^{\otimes
n-1}\hookrightarrow {\cal L}$ (\`{a} homoth\'{e}tie pr\`{e}s) et
${\cal L}/{\cal L}'$ respectivement.

Si $E$ est un syst\`{e}me local sur $X$, on d\'{e}finit le foncteur
$$
M_{n,E}^{\prime d}:D_{{\rm c}}^{{\rm b}}(\mathop{\rm
Fib}\nolimits_{n}^{\prime\, n(n-1)(g-1)}(X))
\rightarrow D_{{\rm c}}^{{\rm b}}(\mathop{\rm
Fib}\nolimits_{i}^{\prime\, n(n-1)(g-1)+d}(X))
$$
par
$$
M_{n,E}^{\prime d}(K)=p_{\ast}p'^{\ast}(K\otimes q^{\ast}{\cal
L}_{E}^{d})[nd]
$$
o\`{u} ${\cal L}_{E}^{d}$ est le faisceau pervers sur $\mathop{\rm
Coh}\nolimits_{0}^{d}(X)$ construit dans la section $3$.

\thm D\'{E}FINITION 4.1
\enonce
Si $E$ est un syst\`{e}me local sur $X$, on d\'{e}finit le complexe
$\mathop{\rm Aut}\nolimits_{E}^{\prime}$ sur $\mathop{\rm
Fib}\nolimits_{n}^{\prime}(X)$ par
$$
\mathop{\rm Aut}\nolimits_{E}^{\prime}|\mathop{\rm Fib}
\nolimits_{n}^{\prime n(n-1)(g-1)+d}(X)=M_{n,E}^{\prime d}(L_{n})
$$
si $d\geq 0$ et
$$
\mathop{\rm Aut}\nolimits_{E}^{\prime}|\mathop{\rm Fib}
\nolimits_{n}^{\prime n(n-1)(g-1)+d}(X)=(0)
$$
si $d<0$.
\endthm
\vskip 5mm

\centerline{\bf 5.  CONSTRUCTION PAR TRANSFORMATIONS DE FOURIER}
\vskip 5mm

Pour chaque $i=1,\ldots ,n$, soient $\mathop{\rm Coh}\nolimits_{i}(X)$
le champ alg\'{e}brique des ${\cal O}_{X}$-Modules coh\'{e}rents
${\cal M}_{i}$ de rang g\'{e}n\'{e}rique $i$ et $\mathop{\rm
Coh}\nolimits_{i}^{\prime}(X)$ le champ alg\'{e}brique des ${\cal
O}_{X}$-Modules coh\'{e}rents ${\cal M}_{i}$ de rang g\'{e}n\'{e}rique
$i$ munis d'une injection ${\cal O}_{X}$-lin\'{e}aire
$s_{i}:(\Omega_{X}^{1})^{\otimes i-1}\hookrightarrow {\cal M}_{i}$
\`{a} homoth\'{e}tie pr\`{e}s.  Ces champs admettent pour ouverts
denses les champs $\mathop{\rm Fib}\nolimits_{i}(X)$ et $\mathop{\rm
Fib}\nolimits_{i}^{\prime}(X)$ et, tout comme ces derniers champs,
leurs composantes connexes $\mathop{\rm Coh}\nolimits_{i}^{d}(X)$ et
$\mathop{\rm Coh}\nolimits_{i}^{\prime\,d}(X)$ sont d\'{e}coup\'{e}es
par le degr\'{e} $d$ de ${\cal M}_{i}$.

\thm LEMME 5.1
\enonce
Il existe une constante $c(g,n)$ qui a la propri\'{e}t\'{e} suivante :
pour tout entier $d\geq c(g,n)$ et tout ${\cal L}\in \mathop{\rm
Fib}\nolimits_{n}^{d}(X)$ tel que $\mathop{\rm Hom}\nolimits_{{\cal
O}_{X}}({\cal L},(\Omega_{X}^{1})^{\otimes n+1})\not =(0)$, ${\cal L}$ est
tr\`{e}s instable au sens o\`{u} il existe une d\'{e}composition non
triviale en somme directe ${\cal L}={\cal L}_{1}\oplus {\cal L}_{2}$
pour laquelle $\mathop{\rm Ext}\nolimits_{{\cal O}_{X}}^{1}({\cal
L}_{1},{\cal L}_{2})=(0)$.
\hfill\hfill$\square$
\endthm

Pour chaque $i=0,\ldots ,n$, soient ${\cal C}_{i}\subset
\mathop{\rm Coh}\nolimits_{i}(X)$ l'ouvert form\'{e} des ${\cal
M}_{i}$ de degr\'{e} $\geq c(n,g)+i(i-1)(g-1)$ tels que $\mathop{\rm
Hom}\nolimits_{{\cal O}_{X}}({\cal M}_{i},(\Omega_{X}^{1})^{\otimes
n+1})=(0)$, et donc a fortiori tels que $\mathop{\rm
Hom}\nolimits_{{\cal O}_{X}}({\cal M}_{i}, (\Omega_{X}^{1})^{\otimes
i})=(0)$ et que $\mathop{\rm Ext}\nolimits_{{\cal
O}_{X}}^{1}((\Omega_{X}^{1})^{\otimes i-1},{\cal M}_{i})=(0)$ par
dualit\'{e} de Serre.  Soient ${\cal V}_{i}$ le champ alg\'{e}brique
des couples $({\cal M}_{i},s_{i})$ o\`{u} ${\cal M}_{i}\in {\cal
C}_{i}$ et $s_{i}\in \mathop{\rm Hom}\nolimits_{{\cal
O}_{X}}((\Omega_{X}^{1})^{\otimes i-1},{\cal M}_{i})$ et ${\cal
V}_{i}^{\vee}$ le champ alg\'{e}brique des extensions
$$
0\rightarrow (\Omega_{X}^{1})^{\otimes i}\rightarrow {\cal M}_{i+1}
\rightarrow {\cal M}_{i}\rightarrow 0
$$
de ${\cal O}_{X}$-Modules coh\'{e}rents avec ${\cal M}_{i}\in {\cal
C}_{i}$.  On a des projections naturelles $\pi_{i}:{\cal V}_{i}
\rightarrow {\cal C}_{i}$ et $\pi_{i}^{\vee}:{\cal V}_{i}^{\vee}
\rightarrow {\cal C}_{i}$ qui sont des fibr\'{e}s vectoriels en
dualit\'{e}.

Introduisons les ouverts
$$
{\cal V}_{i}^{\circ}=\{({\cal M}_{i},s_{i})\mid s_{i}\hbox{ est
injective}\}\subset {\cal V}_{i}
$$
et
$$
{\cal V}_{i}^{\vee\circ}=\{(0\rightarrow (\Omega_{X}^{1})^{\otimes
i}\rightarrow {\cal M}_{i+1} \rightarrow {\cal M}_{i}\rightarrow
0)\mid {\cal M}_{i+1}\in {\cal C}_{i+1}\} \subset {\cal V}_{i}^{\vee}.
$$
Bien entendu, on a un isomorphisme canonique
$$
{\cal V}_{i}^{\vee\circ}\buildrel\sim\over\longrightarrow
{\cal V}_{i+1}^{\circ}
$$
qui envoie $0\rightarrow (\Omega_{X}^{1})^{\otimes i}\rightarrow {\cal
M}_{i+1}\rightarrow {\cal M}_{i}\rightarrow 0$ sur l'injection
$(\Omega_{X}^{1})^{\otimes i}\rightarrow {\cal M}_{i+1}$.

Le groupe multiplicatif agit par homoth\'{e}tie sur les fibr\'{e}s
vectoriels ${\cal V}_{i}$ et ${\cal V}_{i}^{\vee}$ et par passage au
quotient on obtient le diagramme fondamental
$$
\displaylines{\quad\diagram{
&&\overline{{\cal V}}_{n}^{\,\circ}\buildrel\sim\over\leftarrow
\overline{{\cal V}}_{n-1}^{\vee\circ}\hookrightarrow
\overline{{\cal V}}_{n-1}^{\vee}&&&&\overline{{\cal V}}_{n-1}
\hookleftarrow{\overline{\cal V}}_{n-1}^{\,\circ}
\buildrel\sim\over\leftarrow\overline{{\cal V}}_{n-2}^{\vee\circ}
\hookrightarrow&\cr
\noalign{\smallskip}
&\kern -2mm\swarrow\kern -1mm&&\kern -2mm\searrow\kern -2mm
&&\kern -2mm\swarrow\kern -1mm&&\cdots\cr
\noalign{\smallskip}
{\cal C}_{n}&&&&{\cal C}_{n-1}&&&\cr}
\hfill\cr\hfill
\diagram{
&\overline{{\cal V}}_{2}^{\vee}&&&&\overline{{\cal V}}_{2}
\hookleftarrow\overline{{\cal V}}_{2}^{\,\circ}
\buildrel\sim\over\leftarrow\overline{{\cal V}}_{1}^{\vee\circ}
\hookrightarrow\overline{{\cal V}}_{1}^{\vee}&&&&
\overline{{\cal V}}_{1}\hookleftarrow{\overline{\cal V}}_{1}^{\,\circ}
\buildrel\sim\over\leftarrow\overline{{\cal V}}_{0}^{\vee\circ}&&\cr
\noalign{\smallskip}
\cdots &&\kern -2mm\searrow\kern -2mm && \kern -2mm\swarrow\kern
-1mm&&\kern -2mm\searrow\kern -2mm &&\kern -2mm\swarrow\kern
-1mm&&\kern -2mm\searrow\kern -2mm&\cr
\noalign{\smallskip}
&&&{\cal C}_{2}&&&&{\cal C}_{1}&&&&{\cal C}_{0}\cr}\quad}
$$
o\`{u} les fl\`{e}ches obliques sont les projections
$\overline{\pi}_{n}^{\,\circ}:\overline{{\cal V}}_{n}^{\,\circ}=[{\cal
V}_{n}^{\,\circ}/{\Bbb G}_{{\rm m}}]\rightarrow {\cal C}_{n}$,
$\overline{\pi}_{i}:\overline{{\cal V}}_{i}=[{\cal V}_{i}/{\Bbb
G}_{{\rm m}}]\rightarrow {\cal C}_{i}$, $\overline{\pi}_{i}^{\vee}:
\overline{{\cal V}}_{i}^{\vee}=[{\cal V}_{i}^{\vee}/{\Bbb G}_{{\rm
m}}]\rightarrow {\cal C}_{i}$ et $\overline{\pi}_{0}^{\vee\circ}:
\overline{{\cal V}}_{0}^{\vee\circ}= [{\cal V}_{0}^{\vee\circ}/{\Bbb
G}_{{\rm m}}]\rightarrow {\cal C}_{0}$ et o\`{u} les fl\`{e}ches
horizontales sont d'une part les immersions ouvertes
$j_{i}:\overline{{\cal V}}_{i}^{\,\circ}\hookrightarrow
\overline{{\cal V}}_{i}$ et $j_{i}^{\vee}:\overline{{\cal
V}}_{i}^{\vee\circ}\hookrightarrow \overline{{\cal V}}_{i}^{\vee}$
induites par les inclusions ${\cal V}_{i}^{\circ}\subset {\cal V}_{i}$
et ${\cal V}_{i}^{\vee\circ} \subset {\cal V}_{i}^{\vee}$ et d'autre
part les isomorphismes $\iota_{i}:\overline{{\cal V}}_{i}^{\vee\circ}
\buildrel\sim\over \longrightarrow \overline{{\cal
V}}_{i+1}^{\,\circ}$ induits par les isomorphismes ${\cal
V}_{i}^{\vee\circ}\buildrel\sim\over \longrightarrow {\cal
V}_{i+1}^{\circ}$ ci-dessus.

Pour tout syst\`{e}me local $E$ de rang $n$ sur $X$ on d\'{e}finit
alors les complexes $K_{E,i}\in D_{{\rm c}}^{{\rm b}} (\overline{{\cal
V}}_{i}^{\,\circ})$, $i=1,\ldots ,n$, par
$$
K_{E,1}=(\iota_{0})_{\ast}(\overline{\pi}_{0}^{\vee\circ})^{\ast}
{\cal L}_{E},
$$
o\`{u} ${\cal L}_{E}$ est le complexe sur ${\cal C}_{0}$ dont la
restriction \`{a} chaque composante connexe ${\cal C}_{0}^{d}$ de
${\cal C}^{0}$ est le faisceau pervers ${\cal L}_{E}^{d}$ d\'{e}fini
en 3.2, d\'{e}cal\'{e} de $d-1$, et par la relation de r\'{e}currence
$$
K_{E,i+1}=(\iota_{i})_{\ast}(j_{i}^{\vee})^{\ast}\mathop{\rm Four}
\nolimits_{\overline{{\cal V}}_{i}/{\cal C}_{i}}((j_{i})_{!}K_{E,i}),
~i=1,\ldots ,n-1
$$
o\`{u} $\mathop{\rm Four}\nolimits_{\overline{{\cal V}}_{i}/{\cal
C}_{i}}:\mathop{\rm D}_{{\rm c}}^{{\rm b}}(\overline{{\cal
V}}_{i})\rightarrow \mathop{\rm D}_{{\rm c}}^{{\rm
b}}(\overline{{\cal V}}_{i}^{\vee})$ est la transformation de
Fourier homog\`{e}ne pour le fibr\'{e} vectoriel ${\cal
V}_{i}\rightarrow {\cal C}_{i}$ d\'{e}finie dans la section $0$.

Les champs $\overline{{\cal V}}_{i}^{\circ}$ et $\mathop{\rm Fib}
\nolimits_{i}^{\prime}(X)$ sont des ouverts de $\mathop{\rm Coh}
\nolimits_{i}^{\prime}(X)$ et on peut donc consid\'{e}rer leur
intersection.

\thm LEMME 5.2
\enonce
Soit $E$ un syst\`{e}me local de rang $n$ sur $X$.  Pour chaque
$i=1,\ldots ,n$ et chaque entier $d\geq 0$, la restriction de
$K_{E,i}$ \`{a} l'ouvert $\overline{{\cal V}}_{i}^{\circ}\cap
\mathop{\rm Fib} \nolimits_{i}^{\prime\,i(i-1)(g-1)+d}(X)$ de
$\mathop{\rm Coh} \nolimits_{i}^{\prime\,i(i-1)(g-1)+d}(X)$ est
\'{e}gale \`{a}
$$
M_{i,E}^{\prime d}(L_{i})
$$
avec les notations de la section $4$.

En particulier, pour chaque entier $d\geq 0$, les restrictions de
$K_{E,n}$ et de $\mathop{\rm Aut}\nolimits_{E}^{\prime}$ \`{a}
l'ouvert $\overline{{\cal V}}_{n}^{\circ}\cap\mathop{\rm Fib}
\nolimits_{i}^{\prime\,n(n-1)(g-1)+d}(X)$ de $\mathop{\rm Coh}
\nolimits_{i}^{\prime\,n(n-1)(g-1)+d}(X)$ co\"{\i}ncident.
\hfill\hfill$\square$
\endthm
\vskip 5mm

\centerline{\bf 6. LE TH\'{E}OR\`{E}ME D'ANNULATION}
\vskip 5mm

Soient $i\geq 1$ et $d\geq 0$ des entiers.  Consid\'{e}rons le
diagramme
$$\diagram{
\noalign{\vskip 10mm}
\put (16,10){\mathop{\rm Mod}\nolimits_{i}^{d}(X)}
\put (33,10){\smash{\mathop{\hbox to 8mm{\rightarrowfill}}
\limits^{\scriptstyle q}}}
\put (43,10){\mathop{\rm Coh}\nolimits_{0}^{d}(X)}
\arrow(14,8)\dir(-1,-1)\length{8}
\put (8,5){\scriptstyle p'}\cr
\arrow(31,8)\dir(1,-1)\length{8}\cr
\put (36,5){\scriptstyle p}\cr
\kern -1mm \mathop{\rm Fib}\nolimits_{i}(X)&&
\kern 17mm\mathop{\rm Fib}\nolimits_{i}(X)\cr}
$$
o\`{u} $\mathop{\rm Mod}\nolimits_{i}^{d}(X)$ est le champ des couples
$({\cal L},{\cal L}'\subset {\cal L})$ form\'{e} d'un fibr\'{e}
vectoriel de rang $i$ et d'une modification inf\'{e}rieure ${\cal
L}'\subset {\cal L}$ de co-longueur $d$, o\`{u} les projections $p'$
et $p$ envoient $({\cal L},{\cal L}'\subset {\cal L})$ sur ${\cal L}'$
et ${\cal L}$ respectivement et o\`{u} $q$ envoie $({\cal L},{\cal
L}'\subset {\cal L})$ sur le ${\cal O}_{X}$-Module coh\'{e}rent ${\cal
L}/{\cal L}'$ de rang g\'{e}n\'{e}rique $0$ et de longueur $d$.

Pour tout syst\`{e}me local $E$ sur $X$, on a la variante
suivante du foncteur $M_{i,E}^{\prime d}$ de la section 4
$$
M_{i,E}^{d}:D_{{\rm c}}^{{\rm b}}(\mathop{\rm Fib}\nolimits_{i}(X))
\rightarrow D_{{\rm c}}^{{\rm b}}(\mathop{\rm Fib}\nolimits_{i}(X)),~
M_{i,E}^{d}(K)=p_{\ast}p'^{\ast}(K\otimes q^{\ast}{\cal L}_{E}^{d})[id],
$$
o\`{u} ${\cal L}_{E}^{d}$ est le faisceau pervers sur $\mathop{\rm
Coh}\nolimits_{0}^{d}(X)$ d\'{e}fini en $3.2$.

\thm TH\'{E}OR\`{E}ME 6.1 (Gaitsgory, [Ga])
\enonce
Soit $E$ un syst\`{e}me local {\rm irr\'{e}ductible} de rang $n$.
Alors, pour tout $i=1,\ldots ,n-1$ et tout entier $d>in(2g-2)$, le
foncteur $M_{i,E}^{d}$ est identiquement nul.
\hfill\hfill$\square$
\endthm

\rem Remarque $6.2$
\endrem
La d\'{e}monstration de ce th\'{e}or\`{e}me est longue et technique.  Elle
d\'{e}passe largement le cadre de cet expos\'{e}.  Signalons cependant
que, pour $i=1$, cet \'{e}nonc\'{e} n'est autre qu'une reformulation
d'un r\'{e}sultat de Deligne (cf.  [Dr~2]) qui dit que
$$
(\varphi^{(d)})_{\ast}E^{(d)}=0
$$
pour tout syst\`{e}me local irr\'{e}ductible de rang $n>1$ sur $X$ et
tout entier $d>n(2g-2)$, o\`{u} $\varphi^{(d)}:X^{(d)}\rightarrow
\mathop{\rm Pic}\nolimits^{0}(X)$ est le fibr\'{e} projectif
consid\'{e}r\'{e} dans la section 1.
\hfill\hfill$\square$
\vskip 3mm

Frenkel, Gaitsgory et Vilonen utilisent le th\'{e}or\`{e}me 6.1 pour
d\'{e}montrer la proposition suivante.

\thm PROPOSITION 6.3
\enonce
Soit $E$ un syst\`{e}me local {\rm irr\'{e}ductible} de rang $n$.
Alors, pour chaque $i=1,\ldots ,n-1$, la fl\`{e}che canonique
$$
(j_{i})_{!}K_{E,i}\rightarrow (j_{i})_{\ast}K_{E,i}
$$
est un isomorphisme.
\endthm

\rem Id\'{e}e de la d\'{e}monstration
\endrem
Au-dessus de l'ouvert $\mathop{\rm Fib}\nolimits_{i}(X)\cap {\cal
C}_{i}$ de ${\cal C}_{i}$, ${\cal V}_{i}^{\circ}\subset {\cal V}_{i}$
est le compl\'{e}mentaire de la section nulle dans le fibr\'{e} vectoriel
${\cal V}_{i}\rightarrow {\cal C}_{i}$.  D'apr\`{e}s le lemme (0.3),
pour d\'{e}montrer l'assertion au-dessus de cet ouvert, il suffit donc
de montrer que $(\overline{\pi}_{i}^{\,\circ})_{!}K_{E,i}=(0)$ o\`{u}
$\overline{\pi}_{i}^{\,\circ}:\overline{{\cal
V}}_{i}^{\,\circ}\rightarrow {\cal C}_{i}$ est la projection, ou ce
qui revient au m\^{e}me la projection ${\Bbb P}({\cal
V}_{i})\rightarrow {\cal C}_{i}$ puisque les deux projections
co\"{\i}ncident au-dessus de $\mathop{\rm Fib}\nolimits_{i}(X)\cap
{\cal C}_{i}$.  Or, d'apr\`{e}s le lemme 5.2, pour tout entier $d\geq
0$, on a
$$
(\overline{\pi}_{i}^{\,\circ})_{!}K_{E,i}|({\cal C}_{i}\cap \mathop{\rm
Fib}\nolimits_{i}^{i(i-1)(g-1)+d}(X))=M_{i,E}^{d}(\pi_{\ast}'L_{i})|
({\cal C}_{i}\cap \mathop{\rm Fib}\nolimits_{i}^{i(i-1)(g-1)+d}(X))
$$
o\`{u} $\pi':\mathop{\rm Fib}\nolimits_{i}^{\prime}(X)\rightarrow
\mathop{\rm Fib}\nolimits_{i}^{}(X)$ est le morphisme d'oubli de la
section $s:(\Omega_{X}^{1})^{\otimes n-1}\rightarrow {\cal L}$,
d'o\`{u} la conclusion d'apr\`{e}s le
th\'{e}or\`{e}me 6.1.
\hfill\hfill$\square$
\vskip 3mm

Il r\'{e}sulte de la proposition 6.3 que, si $K_{E,i}$ est un faisceau
pervers irr\'{e}ductible, il en est de m\^{e}me de
$(j_{i})_{!}K_{E,i}$ puisque c'est alors le prolongement
interm\'{e}diaire de $K_{E,i}$.  Comme la transformation de Fourier
homog\`{e}ne pr\'{e}serve la perversit\'{e} et
l'irr\'{e}ductibilit\'{e}, on voit par r\'{e}currence sur $i$ que la
restriction de $K_{E,i}$ \`{a} chaque composante connexe de
$\overline{{\cal V}}_{i}^{\,\circ}$ est un faisceau pervers
irr\'{e}ductible.  En particulier, on a d\'{e}montr\'{e}:

\thm COROLLAIRE 6.4
\enonce
Soit $E$ un syst\`{e}me local irr\'{e}ductible de rang $n$.  La
restriction de $K_{E,n}$ \`{a} chaque composante connexe de
$\overline{{\cal V}}_{n}^{\,\circ}$ est un faisceau pervers
irr\'{e}ductible.
\hfill\hfill$\square$
\endthm
\vskip 5mm

\centerline{\bf 7. LA DESCENTE}
\vskip 5mm

Commen\c{c}ons par rappeler deux r\'{e}sultats g\'{e}n\'{e}raux.  Pour
tout sch\'{e}ma (ou champ alg\'{e}brique) $S$ et tout $K\in D_{{\rm
c}}^{{\rm b}}(S)$, notons
$$
\chi_{K}:S\rightarrow {\Bbb Z}
$$
la fonction {\it caract\'{e}ristique d'Euler-Poincar\'{e} de} $K$
d\'{e}finie par
$$
\chi_{K}(s)=\sum_{i}(-1)^{i}\mathop{\rm dim}\nolimits_{C}
H^{i}(K_{s}),~ \forall s\in S,
$$
o\`{u} $K_{s}$ est la fibre de $K$ en $s$.

\thm LEMME 7.1 (Deligne, [Il] Corollaire 2.10)
\enonce
Soit $f:Y\rightarrow Z$ un morphisme propre de sch\'{e}mas, voire un
morphisme repr\'{e}sentable et propre de champs alg\'{e}briques.
Alors, si $K,K'\in D_{{\rm c}}^{{\rm b}}(Y)$ sont localement
isomorphes sur $Y$, les complexes $f_{\ast}K$ et $f_{\ast}K'$ ont
m\^{e}me fonction caract\'{e}ristique d'Euler-Poincar\'{e}
$\chi_{f_{\ast}K}=\chi_{f_{\ast}K'}$ sur $Z$.
\hfill\hfill$\square$
\endthm

\thm LEMME 7.2
\enonce
Soit $S$ un sch\'{e}ma, voire un champ alg\'{e}brique, soit
$V\rightarrow S$ un fibr\'{e} vectoriel de rang constant fini $r$,
d\'{e}finissant un fibr\'{e} projectif $\pi :P={\Bbb P}(V)\rightarrow
S$, et soit $K$ un faisceau pervers irr\'{e}ductible sur $P$.  Alors,
il existe un faisceau pervers {\rm (}n\'{e}cessairement
irr\'{e}ductible{\rm )} $L$ sur $S$ tel que $K$ soit isomorphe \`{a}
$\pi^{\ast}(L)[r-1]$ si et seulement si la fonction
caract\'{e}ristique d'Euler-Poincar\'{e} $\chi_{K}$ est constante le
long des fibres de $\pi$.
\endthm

\rem D\'{e}monstration
\endrem
La partie {\og}{seulement si}{\fg} est triviale.  Supposons donc
$\chi_{K}$ constante le long des fibres de $\pi$.  Par le
th\'{e}or\`{e}me de structure des faisceaux pervers irr\'{e}ductibles,
on a $K\cong i_{\ast}j_{!\ast}E[\mathop{\rm dim}Q]$ o\`{u}
$i:Q\hookrightarrow P$ est le ferm\'{e} irr\'{e}ductible support de
$K$, $j:Q^{\circ}\hookrightarrow Q$ est un ouvert dense et $E$ est un
syst\`{e}me local irr\'{e}ductible sur $Q^{\circ}$.

On v\'{e}rifie dans un premier temps que $Q=\pi^{-1}(T)$ pour un
ferm\'{e} irr\'{e}ductible $T$ de $S$, puis que l'on peut choisir
$Q^{\circ}$ de la forme $\pi^{-1}(T^{\circ})$ pour un ouvert dense
$T^{\circ}$ de $T$ et enfin que $E$ provient d'un syst\`{e}me local
irr\'{e}ductible $F$ sur $T^{\circ}$ puisque les fibres de $\pi$ sont
simplement connexes.

Le faisceau pervers $L$ sur $S$ cherch\'{e} est alors le prolongement
par z\'{e}ro \`{a} $S$ tout entier du prolongement interm\'{e}diaire de
$F[\mathop{\rm dim}T]$ \`{a} $T$.
\hfill\hfill$\square$
\vskip 3mm

\'{E}tudions maintenant la descente de $\mathop{\rm Aut}
\nolimits_{E}^{\prime}$ en un faisceau pervers $\mathop{\rm Aut}
\nolimits_{E}$ sur $\mathop{\rm Fib}\nolimits_{n}(X)$.  Pour cela,
consid\'{e}rons le champ ${\cal Y}_{n}^{d}$ introduit par Drinfeld,
qui classifie les couples $({\cal L},s^{\bullet})$ o\`{u} ${\cal L}$
est un fibr\'{e} vectoriel de degr\'{e} $n(n-1)(g-1)+d$ et
$s^{\bullet}=(s^{1},\ldots ,s^{n})$ est une suite d'injections de
${\cal O}_{X}$-Modules
$$
\matrix{(\Omega_{X}^{1})^{\otimes n-1}\,\smash{\mathop{
\lhook\joinrel\mathrel{\hbox to 8mm{\rightarrowfill}}
}\limits^{\scriptstyle s^{1}}}\,{\cal L}\cr
\cdots\cr
\noalign{\smallskip}
(\Omega_{X}^{1})^{\otimes (n-1)+\cdots +(n-i)}\,\smash{\mathop{
\lhook\joinrel\mathrel{\hbox to 8mm{\rightarrowfill}}
}\limits^{\scriptstyle s^{i}}}\,\bigwedge^{i}{\cal L}\cr
\cdots\cr
\noalign{\smallskip}
(\Omega_{X}^{1})^{\otimes {n(n-1)\over 2}}\,\smash{\mathop{
\lhook\joinrel\mathrel{\hbox to 8mm{\rightarrowfill}}
}\limits^{\scriptstyle s^{n}}}\,\bigwedge^{n}{\cal L}\cr}
$$
qui satisfont les relations de Pl\"{u}cker faisant que la suite
$s^{\bullet}$ d\'{e}finisse un drapeau complet de sous-espaces
vectoriels dans la fibre de ${\cal L}$ au point g\'{e}n\'{e}rique de
$X$.

Le champ ${\cal Y}_{n}^{d}$ est une {\og}{compactification
partielle}{\fg} du champ des triplets $({\cal L},{\cal
L}_{\bullet},s_{\bullet})$ o\`{u}
$$
{\cal L}_{\bullet}=((0)={\cal L}_{0}\subset {\cal L}_{1}\subset\cdots
\subset {\cal L}_{n})
$$
est un drapeau de sous-${\cal O}_{X}$-Modules de ${\cal L}$ et
$s_{\bullet}$ est une suite d'isomorphismes de ${\cal O}_{X}$-Modules
$s_{j}:(\Omega_{X}^{1})^{\otimes j-1}\buildrel\sim\over\longrightarrow
{\cal L}_{n-j+1}/{\cal L}_{n-j}$ pour $j=1,\ldots ,n$.  En effet, le
morphisme
$$
({\cal L},{\cal L}_{\bullet},s_{\bullet})\mapsto ({\cal L},s^{\bullet})
$$
d\'{e}fini par
$$\displaylines{
\quad s^{i}=s_{n}\otimes s_{n-1}\otimes\cdots \otimes
s_{n-i+1}:(\Omega_{X}^{1})^{\otimes (n-1)+(n-2)+\cdots
+(n-i)}
\hfill\cr\hfill
\buildrel\sim\over\longrightarrow {\cal L}_{1}\otimes ({\cal
L}_{2}/{\cal L}_{1})\otimes\cdots \otimes ({\cal L}_{i}/{\cal
L}_{i-1})\cong\bigwedge_{}^{i}{\cal L}_{i}\quad}
$$
identifie ce champ des triplets \`{a} l'ouvert ${\cal Y}_{n}^{d\,
\circ}\subset {\cal Y}_{n}^{d}$ form\'{e} des $({\cal L},
s^{\bullet})$ tels que les conoyaux des injections
$s^{i}:(\Omega_{X}^{1})^{\otimes (n-1)+ \cdots +(n-i)}\hookrightarrow
\bigwedge^{i}{\cal L}$ pour $i=1,\ldots ,n-1$ n'aient pas de torsion.

Le tore ${\Bbb G}_{{\rm m}}^{n}$ agit sur ${\cal Y}_{n}^{d}$
par
$$
(t_{1},\ldots ,t_{n})\cdot ({\cal L},(s^{1},\ldots ,s^{n}))=({\cal
L},(t_{1}s^{1},\ldots ,t_{n}s^{n})).
$$
Cette action respecte l'ouvert ${\cal Y}_{n}^{d\,\circ}$ et est
donn\'{e}e sur le champ des triplets par
$$
(t_{1},\ldots ,t_{n})\cdot ({\cal L},{\cal L}_{\bullet}, (s_{1},\ldots
,s_{n-1},s_{n}))=({\cal L},{\cal L}_{\bullet}, (t_{n}
t_{n-1}^{-1}s_{1},\ldots ,t_{2}t_{1}^{-1}s_{n-1},t_{1}s_{n})).
$$
En particulier, on peut identifier le champ alg\'{e}brique ${\cal
U}_{n}$ (resp.  $\overline{{\cal U}}_{n}=[{\cal U}_{n}/{\Bbb G}_{{\rm
m}}]$) introduit dans la section 4, au quotient de ${\cal
Y}_{n}^{0\,\circ}$ par l'action de ${\Bbb G}_{{\rm m}}$ (resp.  ${\Bbb
G}_{{\rm m}}^{2}$) \`{a} travers le plongement ${\Bbb G}_{{\rm
m}}\hookrightarrow {\Bbb G}_{{\rm m}}^{n},~t\rightarrow
(t,t^{2},\ldots, t^{n})$ (resp.  ${\Bbb G}_{{\rm m}}^{2}
\hookrightarrow {\Bbb G}_{{\rm m}}^{n},~(t,t')\mapsto
(t,t^{2}t',t^{3}t'^{3},\ldots, t^{n}t'^{{n(n-1)\over 2}})$).

\thm PROPOSITION 7.3 ([B-G] Proposition 1.2.2)
\enonce
Le morphisme d'oubli $({\cal L},s^{\bullet})\rightarrow
({\cal L},s^{1})$ passe au quotient en un morphisme de champs
alg\'{e}briques
$$
\widetilde{f}:\widetilde{{\cal Y}}_{n}^{d}:=[{\cal Y}_{n}^{d}/{\Bbb
G}_{{\rm m}}^{n}]\rightarrow \mathop{\rm Fib}\nolimits_{n}^{\prime
\,n(n-1)(g-1)+d}(X)
$$
qui est repr\'{e}sentable et propre.
\hfill\hfill$\square$
\endthm

Consid\'{e}rons le complexe $(\overline{g}_{n})^{\ast}\Psi$ sur
$\overline{{\cal U}}_{n}=[{\cal Y}_{n}^{0\,\circ}/{\Bbb G}_{{\rm
m}}^{2}]$ introduit dans la section 4, et notons $\Phi\in D_{{\rm
c}}^{{\rm b}} (\widetilde{{\cal Y}}_{n}^{0})$ son image directe \`{a}
supports propres par le morphisme compos\'{e} de l'application
quotient $[{\cal Y}_{n}^{0\,\circ}/{\Bbb G}_{{\rm m}}^{2}]\rightarrow
[{\cal Y}_{n}^{0\,\circ}/{\Bbb G}_{{\rm m}}^{n}]$ et de l'inclusion
$[{\cal Y}_{n}^{0\,\circ}/{\Bbb G}_{{\rm m}}^{n}]\subset [{\cal
Y}_{n}^{0}/{\Bbb G}_{{\rm m}}^{n}]=\widetilde{{\cal Y}}_{n}^{d}$.

Pour tout syst\`{e}me local $E$ de rang $n$ sur $X$, on a une variante
du foncteur $M_{n,E}^{d}$
$$
N_{n,E}^{d}:D_{{\rm c}}^{{\rm b}}(\widetilde{{\cal Y}}_{n}^{0})
\rightarrow D_{{\rm c}}^{{\rm b}}(\widetilde{{\cal
Y}}_{n}^{d}),~K\rightarrow p_{\ast}p'^{\ast}(K\otimes q^{\ast}{\cal
L}_{E}^{d})[nd],
$$
d\'{e}finie \`{a} l'aide du diagramme
$$
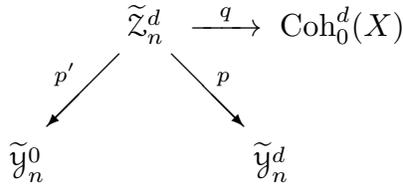
\diagram{
\noalign{\vskip 10mm}
\put (10,10){\widetilde{{\cal Z}}_{n}^{d}}
\put (17,10){\smash{\mathop{\hbox to 8mm{\rightarrowfill}}
\limits^{\scriptstyle q}}}
\put (27,10){\mathop{\rm Coh}\nolimits_{0}^{d}(X)}
\arrow(8,8)\dir(-1,-1)\length{8}
\put (2,5){\scriptstyle p'}\cr
\arrow(15,8)\dir(1,-1)\length{8}\cr
\put (20,5){\scriptstyle p}\cr
\kern -3mm \widetilde{{\cal Y}}_{n}^{0}&&
\kern 16mm \widetilde{{\cal Y}}_{n}^{d}\cr}
$$
o\`{u}
$$
\widetilde{{\cal Z}}_{n}^{d}=\widetilde{{\cal Y}}_{n}^{0}
\times_{\mathop{\rm Fib}\nolimits_{n}(X),p'}{\mathop{\rm Mod}
\nolimits_{n}^{d}(X)}
$$
est le quotient par ${\Bbb G}_{{\rm m}}^{n}$ du champ alg\'{e}brique
des triplets $({\cal L}',s^{\prime\,\bullet},{\cal L'}\subset {\cal L})$ avec
$({\cal L}',s^{\prime\,\bullet})\in {\cal Y}_{n}^{0}$ et $({\cal L'}\subset
{\cal L})\in {\mathop{\rm Mod}\nolimits_{i}^{d}(X)}$, et o\`{u} $p$,
$p'$ et $q$ sont induits par les fl\`{e}ches qui envoient un tel
triplet sur $({\cal L},s^{\bullet})$, $({\cal L}',s^{\prime\,\bullet})$ et
${\cal L}/{\cal L}'$ respectivement avec $s^{i}$ le compos\'{e} de
$s^{\prime\,i}$ et de l'inclusion $\bigwedge^{i}{\cal L}'\subset
\bigwedge^{i}{\cal L}$.

\thm LEMME 7.4
\enonce
Soit $E$ un syst\`{e}me local de rang $n$ sur $X$.  Pour tout entier
$d\geq 0$, on a
$$
\mathop{\rm Aut}\nolimits_{E}^{\prime}|\mathop{\rm Fib}
\nolimits_{n}^{\prime\,n(n-1)(g-1)+d}=\widetilde{f}_{\ast}
N_{n,E}^{d}(\Phi ).
$$
\hfill\hfill$\square$
\endthm

Pour montrer que $\mathop{\rm Aut}\nolimits_{E}^{\prime}$ se descend
par le morphisme $\mathop{\rm Fib}\nolimits_{n}^{\prime}(X)
\rightarrow\mathop{\rm Fib}\nolimits_{n}(X)$ d'oubli de la section,
Frenkel, Gaitsgory et Vilonen utilisent le lemme 7.4, la proposition
7.3 et le lemme 7.1 pour montrer que la fonction caract\'{e}ristique
d'Euler-Poincar\'{e} de $\mathop{\rm Aut}\nolimits_{E}^{\prime}$ ne
d\'{e}pend pas de $E$.  Puis ils calculent cette fonction
d'Euler-Poincar\'{e} quand $E$ est le syst\`{e}me local trivial de
rang $n$, et enfin ils utilisent le lemme 7.2 pour conclure.
\vfill
\eject

\centerline{\bf  BIBLIOGRAPHIE}
\vskip 5mm

\newtoks\ref	\newtoks\auteur
\newtoks\titre	\newtoks\editeur
\newtoks\annee	\newtoks\revue
\newtoks\tome	\newtoks\pages
\newtoks\reste	\newtoks\autre

\def\bibitem#1{\parindent=20pt\itemitem{#1}\parindent=12pt}

\def\livre{\bibitem{[\the\ref]}%
\the\auteur ~-- {\sl\the\titre}, \the\editeur, ({\the\annee}).
\smallskip\smallskip\filbreak}

\def\article{\bibitem{[\the\ref]}%
\the\auteur ~-- \the\titre, {\sl\the\revue} {\the\tome},
({\the\annee}), \the\pages.\smallskip\filbreak}

\def\autre{\bibitem{[\the\ref]}%
\the\auteur ~-- \the\reste.\smallskip\filbreak}

\ref={B-B-D}
\auteur={A. {\pc BEILINSON}, J. {\pc BERNSTEIN}, P. {\pc DELIGNE}}
\reste={Faisceaux pervers, {\it Ast\'{e}risque} {\bf 100}, (1982)}
\autre

\ref={B-D}
\auteur={A. {\pc BEILINSON}, V. {\pc DRINFELD}}
\reste={Quantization of Hitchin's Integrable System and Hecke
Eigensheaves, pr\'{e}publication, http://www.math.uchicago.edu/
$\sim$benzvi/BD/hitchin.ps.gz}
\autre

\ref={B-L}
\auteur={J. {\pc BERNSTEIN}, V. {\pc LUNTS}}
\reste={Equivariant Sheaves and Functors, Lecture Notes in
Mathematics {\bf 1578}, Springer-verlag (1974)}
\autre

\ref={Br}
\auteur={J.-L. {\pc BRYLINSKI}}
\reste={Transformations canoniques, Dualit\'{e} projective,
Th\'{e}orie de Lefschetz, Transformations de Fourier et Sommes
trigonom\'{e}triques, dans {\it G\'{e}om\'{e}trie et Analyse
Microlocales, Ast\'{e}risque} {\bf 140-141}, (1986), 3-134}
\autre

\ref={B-G}
\auteur={A. {\pc BRAVERMAN}, D. {\pc GAITSGORY}}
\reste={Geometric Eisenstein series, pr\'{e}publication,
http://arXiv.org/abs/math/9912097, (1999)}
\autre

\ref={Dr~1}
\auteur={V. G. {\pc DRINFELD}}
\reste={Langlands' conjecture for $\mathop{\rm GL}(2)$ over functional
fields, Proceedings of the International Congress of Mathematicians
(Helsinki, 1978), {\it Acad. Sci. Fennica, Helsinki}, (1980), 565-574}
\autre

\ref={Dr~2}
\auteur={V. G. {\pc DRINFELD}}
\titre={Two-dimensional $l$-adic representations of the fundamental group
of a curve over a finite field and automorphic forms on $\mathop{\rm
GL}(2)$}
\revue={Amer. J. Math.}
\tome={105}
\annee={1983}
\pages={85-114}
\article

\ref={Dr~3}
\auteur={V. G. {\pc DRINFELD}}
\titre={Two-dimensional $l$-adic representations of the Galois group of a
global field of characteristic $p$ and automorphic forms on $\mathop{\rm
GL}(2)$}
\revue={J. of Soviet Math.}
\tome={36}
\annee={1987}
\pages={93-105}
\article

\ref={F-G-K-V}
\auteur={E. {\pc FRENKEL}, D. {\pc GAITSGORY}, D. {\pc KAZHDAN}, K.
{\pc VILONEN}}
\titre={Geometric realization of Whittaker
functions and the Langlands conjecture}
\revue={J. Amer. Math. Soc.}
\tome={11}
\annee={1998}
\pages={451-484}
\article

\ref={F-G-V~1}
\auteur={E. {\pc FRENKEL}, D. {\pc GAITSGORY}, K. {\pc VILONEN}}
\titre={Whittaker patterns in the geometry of moduli spaces of bundles
on curves}
\revue={Annals of Math.}
\tome={153}
\annee={2001}
\pages={699-748}
\article

\ref={F-G-V~2}
\auteur={E. {\pc FRENKEL}, D. {\pc GAITSGORY}, K. {\pc VILONEN}}
\titre={On the geometric Langlands conjecture}
\revue={J. Amer. Math. Soc.}
\tome={15}
\annee={2002}
\pages={367-417}
\article

\ref={Ga}
\auteur={D. {\pc GAITSGORY}}
\reste={On a vanishing conjecture appearing in the geometric Langlands
correspondence, pr\'{e}publication,
http://arXiv.org/abs/math/0204081, (2002)}
\autre

\ref={Il}
\auteur={L. {\pc ILLUSIE}}
\titre={Th\'{e}orie de Brauer et caract\'{e}ristiques
d'Euler-Poincar\'{e} d'apr\`{e}s Deligne}
\revue={Ast\'{e}risque}
\tome={82-83}
\annee={1981}
\pages={161-172}
\article

\ref={La}
\auteur={L. {\pc LAFFORGUE}}
\titre={Chtoucas de Drinfeld et correspondance de Langlands}
\revue={Invent. Math.}
\tome={147}
\annee={2002}
\pages={1-241}
\article

\ref={Lau~1}
\auteur={G. {\pc LAUMON}}
\titre={Correspondance de Langlands g\'{e}om\'{e}trique pour les corps
de fonctions}
\revue={Duke Math. J.}
\tome={54}
\annee={1987}
\pages={309-359}
\article

\ref={Lau~2}
\auteur={G. {\pc LAUMON}}
\reste={Faisceaux automorphes pour $\mathop{\rm GL}\nolimits_{n}$: la
premi\`{e}re construction de Drinfeld, pr\'{e}publication,
http://arXiv.org/abs/math/9511004, (1995)}
\autre

\ref={Lau~3}
\auteur={G. {\pc LAUMON}}
\reste={Transformation de Fourier homog\`{e}ne, en pr\'{e}paration}
\autre

\ref={L-M}
\auteur={G. {\pc LAUMON}, L. {\pc MORET}-{\pc BAILLY}}
\titre={Champs alg\'{e}briques}
\editeur={Springer-Verlag}
\annee={1999}
\livre

\ref={Ly~1}
\auteur={S. {\pc LYSENKO}}
\titre={Local geometrized Rankin-Selberg method for $\mathop{\rm GL}(n)$}
\revue={Duke Math. J.}
\tome={111}
\annee={2002}
\pages={451-493}
\article

\ref={Ly~2}
\auteur={S. {\pc LYSENKO}}
\reste={Global geometrized Rankin-Selberg method for $\mathop{\rm GL}(n)$,
http:// arxiv.org/abs/math.AG/0108208, (2001)}
\autre

\ref={Se}
\auteur={J.-P. {\pc SERRE}}
\titre={Groupes alg\'{e}briques et corps de classes}
\editeur={Hermann}
\annee={1975}
\livre
\vskip 7mm

\line{\hfill\hfill\hbox{{\vtop{\tabalign G\'{e}rard Laumon\cr
\tabalign Universit\'{e} Paris-Sud\cr
\tabalign et CNRS, UMR 8628\cr
\tabalign Math\'ematique, B\^{a}t. 425\cr
\tabalign F-91405 Orsay Cedex (France)\cr
\tabalign Gerard.Laumon@math.u-psud.fr\cr}}}}
\bye